\documentclass[11pt, a4paper]{article}
\usepackage[fleqn]{amsmath} 
\usepackage{amssymb} 
\usepackage{fontenc}
\usepackage[dvips]{graphicx}
\usepackage{epsfig}
\usepackage{booktabs}
\usepackage{undertilde}
\usepackage{color}
\usepackage{geometry} 
\usepackage{fancyhdr}
\usepackage{flafter}

\begin{document}
\vspace{0.5cm}
\begin{center} {\bfseries\Large Coarse Variables of Autonomous ODE Systems and Their Evolution}\\
\vspace{0.5cm}
{\bfseries\large Likun Tan, Amit Acharya\footnote{Tel.: +1 412 268 4566; fax: +1 412 268 7813. \par
                                             \emph{E-mail address}: acharyaamit@cmu.edu}, Kaushik Dayal}\\
\vspace{0.3cm}
{\large Civil and Environmental Engineering, Carnegie Mellon University}\\
\end{center}
\vspace{0.5cm}

\thispagestyle{empty}
\noindent
{\bfseries Abstract.} Given an autonomous system of ordinary differential equations (ODE), we consider developing practical models for the deterministic, slow/coarse behavior of the ODE system. Two types of coarse variables are considered. The first type consists of running finite time averages of phase functions. Approaches to construct the coarse evolution equation for this type are discussed and implemented on a `Forced' Lorenz system and a singularly perturbed system whose fast flow does not necessarily converge to an equilibrium. We explore two strategies. In one, we compute (locally) invariant manifolds of the fast dynamics, parameterized by the slow variables. In the other, the choice of our coarse variables automatically guarantees them to be `slow' in a precise sense. This allows their evolution to be phrased in terms of averaging utilizing limit measures (probability distributions) of the fast flow. Coarse evolution equations are constructed based on these approaches and tested against coarse response of the `microscopic' models. The second type of coarse variables are defined as (non-trivial) scalar state functions that are required by design to evolve autonomously, to the extent possible, with the goal of being candidate state functions for unambiguously initializable coarse dynamics. The question motivates a mathematical restatement in terms of a first-order PDE. A computational approximation is developed and tested on the Lorenz system and the Hald Hamiltonian system.  \\

\noindent
\emph{Keywords:} Coarse-graining; Time-averaging; ODE; Multi-scale modeling\\


\noindent
{\bfseries 1. Introduction} \\

\noindent
Obtaining coarse response of a system of ODEs containing rapidly oscillatory or decaying (stiff) response without detailed information on the evolution of the original (fine) variables is an interesting, but challenging task. The coarse response of the fine dynamics is sought with the following (ideal) objectives: (1) For oscillatory fine response, the fundamental frequency of the coarse dynamics should be orders of magnitude less than that of the fine dynamics. (2) The coarse dynamics should be autonomous in the sense that its evolution is independent of fine initial conditions. One possible application of our work would be in Molecular Dynamics (MD) simulations. MD is based on a particle description of atoms or molecules and treats the interaction between the particles using Newtonian dynamics. However, a key difficulty in MD is the extremely large separation between the time-scales of atomic bond vibrations (femtoseconds) and even the smallest time-scales of engineering interest (nano to micro seconds). This limits MD to unrealistically small time steps for understanding coarse evolution at engineering time-scales. In multi-scale analysis, it is beneficial if in the time integration of the system, the time step in the slow/coarse dynamics is much larger than in the atomistic dynamics. The objective of this work is to develop ideas towards the computational coarse-graining of such dynamical systems. In particular, we focus on methods for defining good `slow' variables for probing macroscopic response and defining their closed evolution (i.e. in terms of themselves), to the extent possible. \par

The first part of our work relates to defining the evolution of coarse variables. These objects are the `slow' state variables in a system with an explicit split into slow and fast motions (e.g. constants of motion for the fast flow), and/or time averaged functions of fast variables. The latter have the interesting property that their evolution is necessarily slow on the fast time scale, as shown in~\cite{acharya2010coarse}. Thus, whether a-priori available, or induced by time-averaging, we deal with singularly perturbed systems of autonomous ODE. The approaches that we apply to define the coarse evolution equation are the following:
\begin{itemize}
\item The first approach is to generate an augmented fine system with the introduction of time-delayed fine variables (where time delay is the averaging period), and compute an appropriate collection of local/overflowing invariant manifolds of the discrete-delay augmented fine ODE system as graphs parameterized by the coarse variables a priori. With the computed collection of invariant manifolds, a closed dynamics for the coarse variables can be defined. We refer to this procedure as Parameterized Locally Invariant Manifolds (PLIM)~\cite{acharya2005parametrized, acharya2006computational, sawant2006model, acharya2010coarse}. The invariant manifolds (computed locally) are considered as a representation of the fine variables on the coarse phase space.
\item The second question we explore originates from the reduced order approach (the Tikhonov approach~\cite{tikhonov1952systems}) that utilizes the differential-algebraic equations (DAE) representing slow manifold(s) to approximate the limit dynamics. The idea is to check whether, even when the fast flow does not converge to an equilibrium, the slow evolution of the first moment of the invariant measure generated by the fast flow for fixed slow variables is well approximated by the DAE of the Tikhonov method.
\item The third approach may be viewed as a natural extension of the method of averaging motivated by practical considerations, naturally for systems whose limiting fast flows may not converge to an equilibrium. The theory of~\cite{artstein1996singularly, artstein2007slow} is applied here, where the limit behavior of the fast flow (parameterized by the slow variable) in a singulary perturbed system is characterized by an invariant measure rather than algebraic equations which are used in the Tikhonov approach. However, this theory requires that variables that can be averaged be `slow', i.e. orthogonal to the fast flow in a precise sense, but does not provide a prescription for generating such `orthogonal observables'. Following~\cite{slemrod2010time} we develop and demonstrate a natural class of slow variables based on time averaged coarse observables, and generate their evolution, where this approach is referred as `Practical Time Averaging'.
\end{itemize}

The second part of this paper relates to the definition and computation of instantaneous (as opposed to time-averaged) coarse variables that aims to allow for an autonomous coarse response. Defining unambiguously initializable coarse dynamics with information only on the coarse state is the ideal goal (i.e. the time derivative of the coarse variable should be uniquely determined given a coarse state). However, this requirement for generating an autonomous coarse response can hardly be satisfied, as noted in the implementation of PLIM (for example, in the work of~\cite{sawant2006model, acharya2006computational}) and the theory of~\cite{artstein1996singularly, artstein2007slow}. In the case of PLIM, the invariant manifolds computed on the coarse phase space are not unique. A physical interpretation is that there could possibly exist multiple fine states that are associated with one single coarse state and these fine states affect the coarse evolution out of the attained coarse state. Thus, one requires knowledge of fine initial conditions to ensure a correct coarse response. In Practical Time Averaging, the `right-hand side' term of the coarse evolution equation is interpreted with an invariant measure which generally cannot be determined by the coarse state only. This second part of the paper develops further the preliminary ideas set forth in~\cite{acharya2007choice} in an attempt to define a class of coarse variables whose dynamical response can be uniquely defined by a given coarse state. \par
This paper is organized as follows: In Section 2, we provide the theoretical framework of PLIM and present the numerical results from the implementation of PLIM on a `Forced' Lorenz system. In Section 3, we show the procedures to write down the coarse evolution equation for singularly perturbed systems based on the Tikhonov approach and the Young measure theory (i.e.~\cite{artstein1996singularly, artstein2007slow, slemrod2010time}). Each of these approaches is followed by numerical studies of model problems. In Section 4, we formulate a first-order PDE to generate a class of coarse variables that formally allow autonomous coarse response. We obtain such coarse variables from solving the PDE and discuss some desirable properties of these coarse variables. We end in Section 5 with some discussions of the strengths and drawbacks of these multi-time scale modeling techniques.\\

\noindent
{\bfseries 2. Parameterized Locally Invariant Manifolds} \\

\noindent
 The main ingredients of PLIM are discussed in~\cite{acharya2005parametrized, acharya2006computational, sawant2006model, acharya2010coarse}. PLIM is based on more-or-less standard ideas for computing invariant manifolds using the invariance equation of an autonomous ODE system~\cite{sacker1965new, carr1983applicationI, carr1983applicationII, muncaster1984invariant, foias1988inertial, gorbana2004constructive} with \emph{one essential, and important, difference}. The coarse `retained' master variables whose dynamics is sought are not some of the fine variables, or a function of those. Instead, they are running time averages of functions defined on fine trajectories. With a simple augmentation of the fine dynamics as proposed in~\cite{acharya2006computational}, a coarse dynamics can be obtained for the time averages. The advantage of this strategy lies in the fact that time-averaging provides a natural way to separate time-scales in the system consisting of the augmented fine ODEs and the coarse variables, even though the original (non-augmented) fine ODE system may not come equipped with an in-built and identifiable separation of scales in its variables~\cite{acharya2010coarse}.\\

\noindent
2.1 THEORETICAL FRAMEWORK \\

\noindent
We are interested in nonlinear, autonomous fine dynamics defined as
\begin{equation}\label{eq:org}
\begin{matrix}
\begin{aligned}
& \frac{\displaystyle df}{\displaystyle dt}(t)=H(f(t),I(t))\quad;\quad \frac{\displaystyle dI}{\displaystyle dt}(t)=\frac{\displaystyle 1}{\displaystyle \mbox{T}}L(I(t))\medskip\\
& f(0)=f_{*}\quad;\quad I(0)=I_{*},
\end{aligned}
\end{matrix}
\end{equation}
where $f$ is an $N$-dimensional vector of fine degrees of freedom and $H$ is a generally nonlinear function of fine states, referred to as the vector field of the fine dynamical system. $N$ can be large in principle, and the function $H$ rapidly oscillating. $I$ is an $n$ dimensional vector of slow loading variables that do not vary appreciably compared with the oscillation of fine variables. Here, T is a non-dimensional parameter that represents the ratio of natural frequency of the fine dynamics to the loading frequency, i.e. $\mbox{T}:=\omega_f/\omega$, where $\omega_f$ denotes a characteristic frequency of the fast dynamics and $\omega$ is the time scale of the loading (e.g. constant rate of monotonic loading, frequency of cyclic loading). $L$ takes the form $L(I):=\hat{L}(I)\omega_f$, where $\hat{L}$ is a bounded smooth function with dimensions of the load vector $I$. For $\omega_f$ fixed (which is what we have in mind), by definition, $\mbox{T}\rightarrow\infty$ as $\omega\rightarrow0$. \par
Let $c$ denote the coarse variables of dimension $m$ and $\Lambda$ be a user-specified function of the fine states producing vectors with $m$ components whose time averages over intervals of period $\tau$ can be measured in principle and are of physical interest. Given the fixed time interval $\tau$ characterizing the resolution of coarse measurements in time, a coarse trajectory corresponding to each fine trajectory $f(\cdot)$ is defined as the following running time average:
\begin{equation}\label{eq:def_cos}
c(t):=\frac{\displaystyle 1}{\displaystyle \tau}\int_t^{t+\tau}\Lambda(f(p))dp,
\end{equation}
where $\tau<1/\omega$. Thus, $\tau$ is the period of time averaging, referred to the fast time scale. Clearly,
\begin{equation}
\frac{\displaystyle dc}{\displaystyle dt}(t)=\frac{\displaystyle 1}{\displaystyle \tau}[\Lambda(f(t+\tau))-\Lambda(f(t))].
\end{equation}
From the definition of time-averaging, forward trajectories $f_f(\cdot)$ and $I_f(\cdot)$ are introduced as
\begin{equation}
\begin{matrix}
\begin{aligned}
& f_{f}(t):=f(t+\tau)\\
& I_{f}(t):=I(t+\tau),
\end{aligned}
\end{matrix}
\end{equation}
then
\begin{equation}
\begin{matrix}
\begin{aligned}
& \frac{\displaystyle df_{f}}{\displaystyle dt}(t)=H(f_{f}(t),I_{f}(t)) \\
& \frac{\displaystyle
dI_{f}}{\displaystyle dt}(t)=\frac{\displaystyle 1}{\displaystyle
\mbox{T}}L(I_{f}(t)).
\end{aligned}
\end{matrix}
\end{equation}
With these definitions in hand, an augmented system is introduced:
\begin{equation}\label{eq:aug}
\begin{matrix}
\begin{aligned}
& \frac{\displaystyle df_{f}}{\displaystyle
dt}(t)=H(f_{f}(t),I_{f}(t))\quad;\quad f_{f}(0)=f_{**}\medskip\\
& \frac{\displaystyle df}{\displaystyle
dt}(t)=H(f(t),I(t))\quad;\quad
f(0)=f_{*}\medskip\\
& \frac{\displaystyle dI_{f}}{\displaystyle
dt}(t)=\frac{\displaystyle 1}{\displaystyle \mbox{T}}L(I_{f}(t))\quad;\quad
I_{f}(0)=I_{**}\medskip\\
& \frac{\displaystyle dI}{\displaystyle dt}(t)=\frac{\displaystyle
1}{\displaystyle \mbox{T}}L(I(t))\quad;\quad I(0)=I_{*}\medskip\\
& \frac{\displaystyle dc}{\displaystyle dt}(t)=\frac{\displaystyle
1}{\displaystyle \tau}[\Lambda(f_{f}(t))-\Lambda(f(t))]\quad;\quad
c(0)=c_{*},
\end{aligned}
\end{matrix}
\end{equation}
where we assume that $f_{*}$ and $f_{**}$ are arbitrarily specifiable. Thus, the set of evolutions described by~(\ref{eq:aug}) is strictly larger than that defined by~(\ref{eq:org}), the latter being a subset of the former.\par

Now we seek functions $G^i_f(c,I,I_f)$ and $G^i(c,I,I_f)$, $i=1,\dots,N$ that satisfy the first order, quasilinear partial differential equations
\begin{equation}\label{eq:gov}
\begin{matrix}
\left. \begin{matrix}
\begin{aligned}
& {\sum_{k=1}^{m}}\frac{
\partial G^J_f}{
\partial c^k} \left(\frac{1}{
\tau}[\Lambda^{k}(G_{f})-\Lambda^{k}(G)] \right)+
{\sum_{i=1}^{n}}\frac{1}{\mbox{T}}\left(\frac{
\partial G^J_f}{\partial I^i}L^{i}(I)+\frac{\partial G^J_f}{\partial
I^i_f}L^{i}(I_f)\right)=H^{J}(G_f,I_f)\medskip \\
& {\sum_{k=1}^{m}}\frac{\partial G^J}{\textstyle
\partial c^k}\left(\frac{1}{\tau}[\Lambda^{k}(G_{f})-\Lambda^{k}(G)]\right)+
{\sum_{i=1}^{n}}\frac{1}{\mbox{T}}\left(\frac{\partial G^J}{\partial
I^i}L^{i}(I)+\frac{\partial G^J}{\partial
I^i_f}L^{i}(I_f)\right)=H^{J}(G,I)
\end{aligned}
\end{matrix} \right \}\medskip \\
\quad {J \thinspace = \thinspace 1\thinspace to\thinspace N,}
\end{matrix}
\end{equation}
at least locally in $c$-space. These are generated by assuming $G_f=f_f$ and $G=f$, and substituting in (\ref{eq:aug}). Suppose we have such a pair of functions over the coarse domain containing the point
\begin{equation}
c_*:=c(0)
\end{equation}
defined from (\ref{eq:org}) and (\ref{eq:def_cos}) which, moreover, satisfies the conditions
\begin{equation}\label{eq:constrain}
\begin{aligned}
& G_f(c_*)=f_{**} \\
& G(c_*)=f_*.
\end{aligned}
\end{equation}
Also, given an initial state ($f_*$, $I_*$) we denote ($f_{**}$, $I_{**}$) as the solution of (\ref{eq:org}) evaluated at time $\tau$. It is easy to see that a local-in-time fine trajectory defined by
\begin{equation}
\begin{aligned}
&\Gamma_f(t):=G_f(c(t))\\
&\Gamma(t):=G(c(t))
\end{aligned}
\end{equation}
through the coarse local trajectory satisfying
\begin{equation}\label{eq:coarse}
\begin{matrix}
\begin{aligned}
& \frac{\displaystyle dc}{\displaystyle dt}=\frac{\displaystyle
1}{\displaystyle \tau}\left[\Lambda\left(G_{f}(c,I,I_{f})\right)-\Lambda\left(G(c,I,I_{f})\right)\right]\medskip\\
& {+\;\mbox{the evolutionary equations for}\; I,\thinspace I_{f}}
\end{aligned}
\end{matrix}
\end{equation}
is the solution of~(\ref{eq:aug}) (locally). A solution pair $(G_f,G)$ of (\ref{eq:gov}) represents a parametrization of a locally invariant manifold of the dynamics~(\ref{eq:aug}). Thus,~(\ref{eq:coarse}) is the closed coarse evolution corresponding to the augmented fine system~(\ref{eq:aug}) and~(\ref{eq:constrain}), \emph{but only locally}. \\

\clearpage
\noindent
2.2 A MODEL PROBLEM: `MONOTONICALLY FORCED' LORENZ SYSTEM \\

\noindent
We test the idea of PLIM on the Lorenz system~\cite{lorenz1963deterministic}, and take the time average of state variables in the fine phase space as coarse variables. As is known, the Lorenz system is a three-dimensional chaotic system evolving over time in a complex, non-repeating pattern within a bounded region of phase space (for appropriate parameter values). The ODE is given as:
\begin{equation}\label{eq:lorenz}
\begin{aligned}
&\dot{x}=\sigma(y-x) \\
&\dot{y}=\gamma x-y-xz \\
&\dot{z}=xy-\beta z \\
&\sigma=10,\;\beta=8/3,\;\gamma=25.
\end{aligned}
\end{equation}
We are interested in the evolution of running time averages of $x$ over different periods $\tau$, while the Lorenz system contains high frequency non-periodic bounded oscillations. These oscillations more or less average out to zero. Thus, for the purpose of testing our ideas, we construct a `Forced' Lorenz system whose fine variables show high frequency oscillations superposed on a monotonically increasing/decreasing profile.\par

Let $\tilde{x}(t)$ denote the trajectory whose time average has a monotonically increasing feature while keeping the high frequency oscillations at the fine scale. Let $L(t)$ be a trajectory that is monotonically increasing. One possible construction of $\tilde{x}$ would be:
\begin{equation}\label{eq:defx}
\tilde{x}(t)=x(t)+L(t), \quad\mbox{where} \quad \dot{L}=\frac{1}{1+L}.
\end{equation}
The trajectory with a general increasing pattern can be obtained by simply considering a linear combination of the original trajectory with a monotonically increasing curve. However, to derive the corresponding ODE requires some more work. From ~(\ref{eq:defx}), we have
\begin{equation}\label{eq:reversex}
x(t)=\tilde{x}(t)-L(t).
\end{equation}
We also have
\begin{equation}\label{eq:forcedlorenzx}
\dot{\tilde{x}}=\dot{x}+\dot{L}=\sigma(y-x)+\frac{1}{1+L}.
\end{equation}
Substituting~(\ref{eq:reversex}) into~(\ref{eq:forcedlorenzx}) and~(\ref{eq:lorenz})$_{2,3}$, we obtain the ODE of the `Forced' Lorenz system:
\begin{equation}
\begin{aligned}
&\dot{\tilde{x}}=\sigma(y-\tilde{x})+\sigma L+\frac{1}{1+L} \\
&\dot{y}=\gamma\tilde{x}-y-\tilde{x}z-\gamma L+Lz \\
&\dot{z}=\tilde{x}y-\beta z-Ly \\
&\dot{L}=\frac{1}{1+L}.
\end{aligned}
\end{equation}
For simplicity, we denote $\tilde{x}$ as $x$ and rewrite the ODE as the `Forced' Lorenz system:
\begin{equation}\label{eq:forcedlorenz}
\begin{aligned}
&\dot{x}=\sigma(y-x)+\sigma L+\frac{1}{1+L} \\
&\dot{y}=\gamma{x}-y-xz-\gamma{L}+Lz \\
&\dot{z}=xy-\beta{z}-Ly  \\
&\dot{L}=\frac{1}{1+L}.
\end{aligned}
\end{equation}
Figure~\ref{Lorenz_fig1} shows the trajectories of the original Lorenz system and the `Forced' Lorenz system in the $(x,y,z)$ phase space. The initial conditions are set to be the same. For the Lorenz system, the trajectories evolve around the steady state points with no preference of moving towards any direction, while for the `Forced' Lorenz system, the trajectories move towards the positive $x$ direction over time, keeping the same amplitude in the $y$ and $z$ directions. \par

\begin{center}
\includegraphics [width=.5\columnwidth]{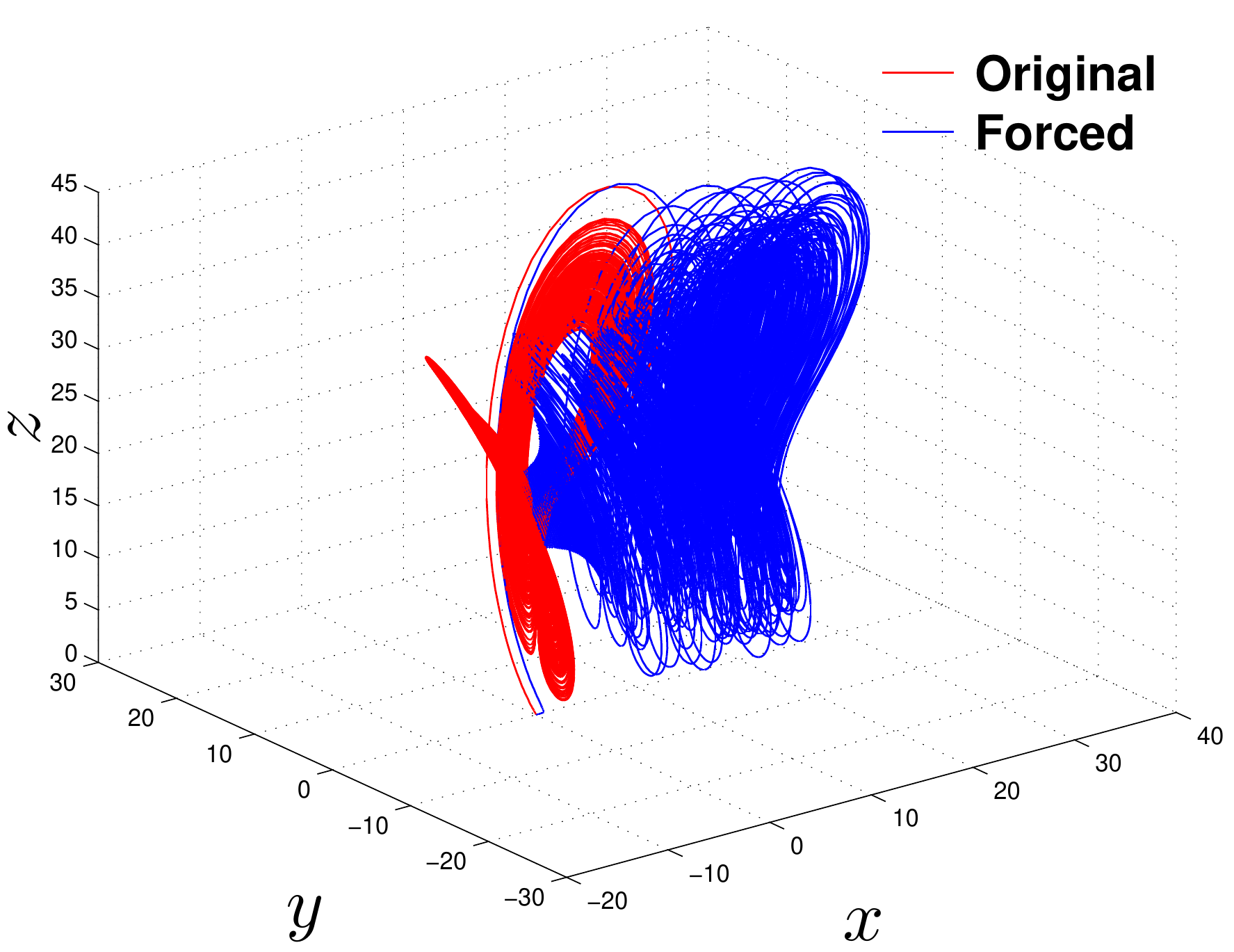}
\makeatletter
\def\@captype{figure}
\makeatother
\caption{Comparison of original Lorenz system and `Forced' Lorenz system}\label{Lorenz_fig1}
\end{center}

\noindent
We take the running time average of $x$ as a coarse variable:
\begin{equation}
\overline{x}(t)=\frac{1}{\tau}\int^{t+\tau}_{t}x(p)dp \Longrightarrow \dot{\overline{x}}(t)=\frac{1}{\tau}\left[x(t+\tau)-x(t)\right]=\frac{1}{\tau}\left[x_f-x\right],\;\;\mbox{where}\;\; x_f:=x(t+\tau).
\end{equation}
We want $\tau$ to be large enough. The highly oscillatory curves can be reduced to a constant by taking an infinitely long time average. Practically, the time averaged variable will be gradually smoothed by increasing the time average, as shown in Figure~\ref{Lorenz_fig2}. The estimated lowest fundamental period of the Lorenz system is around 0.0455 (This approximation is based on calculating the eigenvalues of the linearized system of~(\ref{eq:lorenz}). Note that time is dimensionless in the original definition of Lorenz system~\cite{lorenz1963deterministic}, so this estimated number is also dimensionless). The time step for evolving the fine dynamics is of the order of $10^{-3}$. We set the averaging period to be 50, which is $5*10^4$ times larger than the fine time step.
\begin{center}
\includegraphics [width=.4\columnwidth]{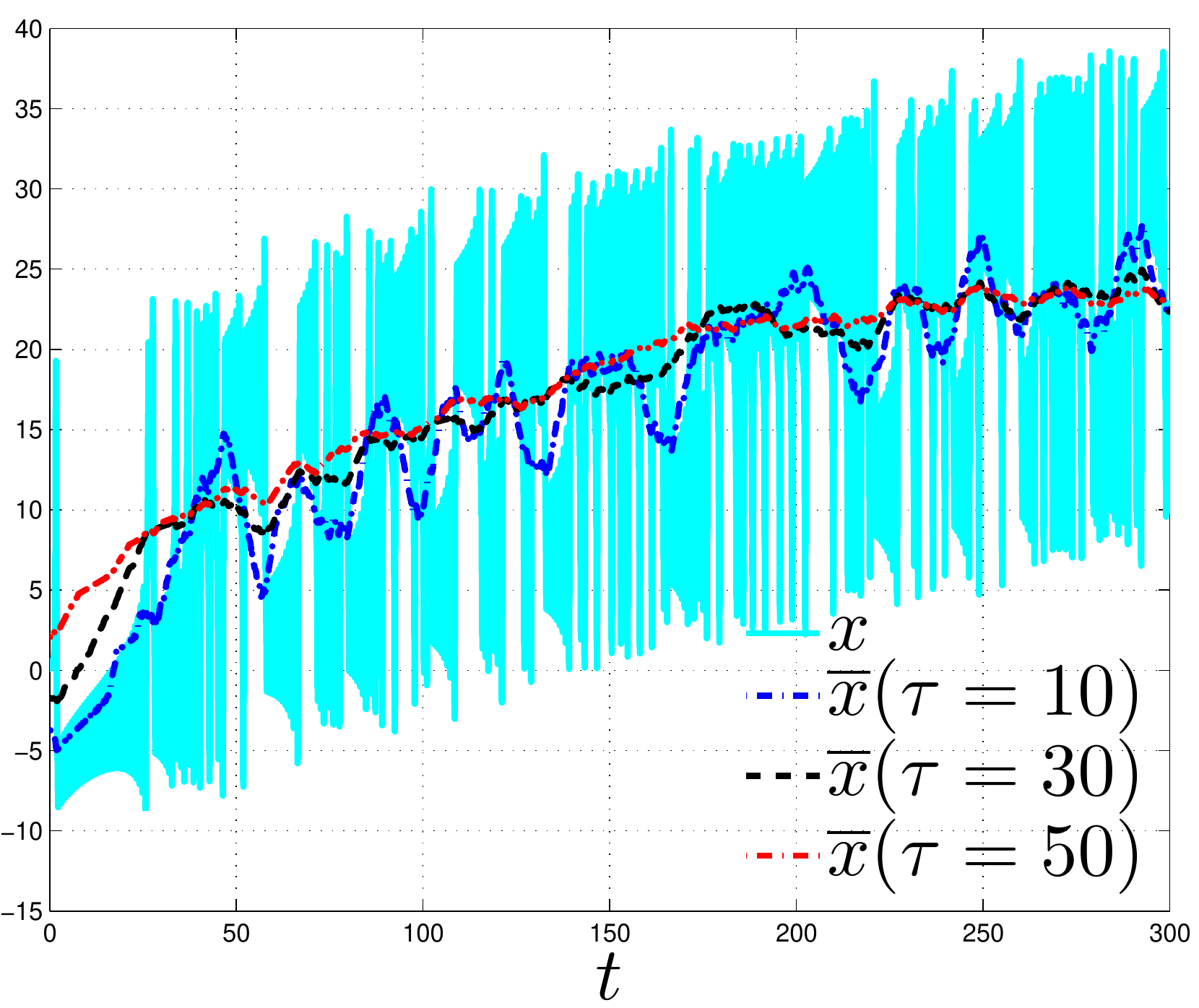}
\makeatletter
\def\@captype{figure}
\makeatother
\caption{Comparison of $x(t)$ and $\overline{x}(t)$ with different time averages}\label{Lorenz_fig2}
\end{center}

\noindent
With the introduction of the time shifted variables $x_f$, $y_f$, $z_f$ and $L_f$, the corresponding augmented system takes the form
\begin{equation}\label{eq:auglorenz}
\begin{aligned}
&\dot{x_f}=\sigma(y_f-x_f)+\sigma{L_f}+\frac{1}{1+L_f} \\
&\dot{y_f}=\gamma{x_f}-y_f-{x_f}{z_f}-\gamma{L_f}+{L_f}{z_f} \\
&\dot{z_f}={x_f}{y_f}-\beta{z_f}-{L_f}{y_f}\\
&\dot{L_f}=\frac{1}{1+L_f}.
\end{aligned}
\end{equation}
Taking $\overline{x}$ and $L$ as the coarse variables, we first construct a mapping on the coarse domain that takes values in the fine phase space, i.e. $x=G_1(\overline{x},L)$, $y=G_2(\overline{x},L)$, $z=G_3(\overline{x},L)$, $x_f=G_4(\overline{x},L)$, $y_f=G_5(\overline{x},L)$, $z_f=G_6(\overline{x},L)$ and $L_f=G_7(\overline{x},L)$. Thus, the governing equations for solving the invariant manifolds have the form
\begin{equation}
\frac{\partial G_i}{\partial\overline{x}}\left(\frac{G_4-G_1}{\tau}\right)+\frac{\partial G_i}{\partial L}\left(\frac{1}{1+L}\right)=H_i(\utilde{G})\quad\forall\;i\;\mbox{from}\:1\: \mbox{to}\: 7,
\end{equation}
where the right-hand side terms are
\begin{equation}
\begin{aligned}
&H_i(\utilde{G})=\sigma(G_{i+1}-G_i)+\sigma{G_{7+(i-1)/3}}+\frac{1}{1+G_{7+(i-1)/3}};\quad i=1\thinspace \mbox{or}\thinspace4 \\
&H_i(\utilde{G})=\gamma{G_{i-1}}-G_i-G_{i-1}G_{i+1}-\gamma{G_{7+(i-2)/3}}+G_{i+1}G_{7+(i-2)/3};\quad i=2\thinspace \mbox{or}\thinspace5 \\
&H_i(\utilde{G})=G_{i-2}G_{i-1}-\beta{G_i}-G_{i-1}G_{7+(i-3)/3};\quad i=3\thinspace \mbox{or}\thinspace6 \\
&H_i(\utilde{G})=\frac{1}{1+G_i}; \quad i=7.
\end{aligned}
\end{equation}
The coarse theory is defined as
\begin{equation}\label{eq:coarse_theory}
\begin{aligned}
&\dot{\overline{x}}=\frac{G_4(\overline{x},L)-G_1(\overline{x},L)}{\tau}\\
&\dot{L}=\frac{1}{1+L}.
\end{aligned}
\end{equation}

\noindent
To obtain the time average of $x$ in the original Lorenz system, we simply subtract the function $L$ from $\overline{x}$ computed from~(\ref{eq:coarse_theory}).\par
Figure~\ref{Lorenz_fig3}(a) shows the comparison between $\overline{x}(t)$ from evolving the fine theory and the coarse theory with different coarse-to-fine time-step ratios (denoted by $c/f$). The coarse trajectories exhibit consistency of a general growing trend, although the two curves do not perfectly agree. When the coarse responses with different time-step ratios are compared, the results agree well up to $c/f=10^4$. This observation is consistent with the previous choice of time average, which is of the order of $10^4$ of the fine scale. Figure~\ref{Lorenz_fig3}(b) is a magnified version of the results from $t$=0 to 10, which demonstrates that the responses are much more consistent at the early steps of evolution. Figure~\ref{Lorenz_fig3}(c) shows a collection of coarse trajectories computed by evolving the fine theory, where the difference of any two coarse initial conditions defined in the $l^2$-norm is smaller than $10^{-3}$, while the difference of the associated fine initial conditions could be quite large. At the coarse scale, the differences among the coarse initial values are considered trivial, nevertheless the trajectories do not totally overlap due to the chaotic nature of the Lorenz system. On the other hand, all of the trajectories follow a gradual increasing pattern as the way it is built up for the `Forced' Lorenz system. The coarse responses from evolving the coarse theory show a better consistency, as seen in Figure~\ref{Lorenz_fig3}(d), where the trajectories are almost overlapping. Here the initial conditions are the same as before. This phenomenon implies \emph{coarse stability} of the constructed coarse model via time-averaging in that small differences of coarse initial values lead to small differences in coarse evolution. The chaotic feature observed at the fine scale is no longer maintained at the coarse scale.\par
Figure~\ref{Lorenz_fig4} shows a group of coarse trajectories from evolving the fine theory, where the initial data are arbitrarily chosen, that is to say, the fine initial conditions corresponding to each coarse trajectory can be significantly different. However, due to the long time averaging, the coarse responses exhibit a similar increasing pattern as noticed in previous results. \par
\begin{center}
\includegraphics [width=.9\columnwidth]{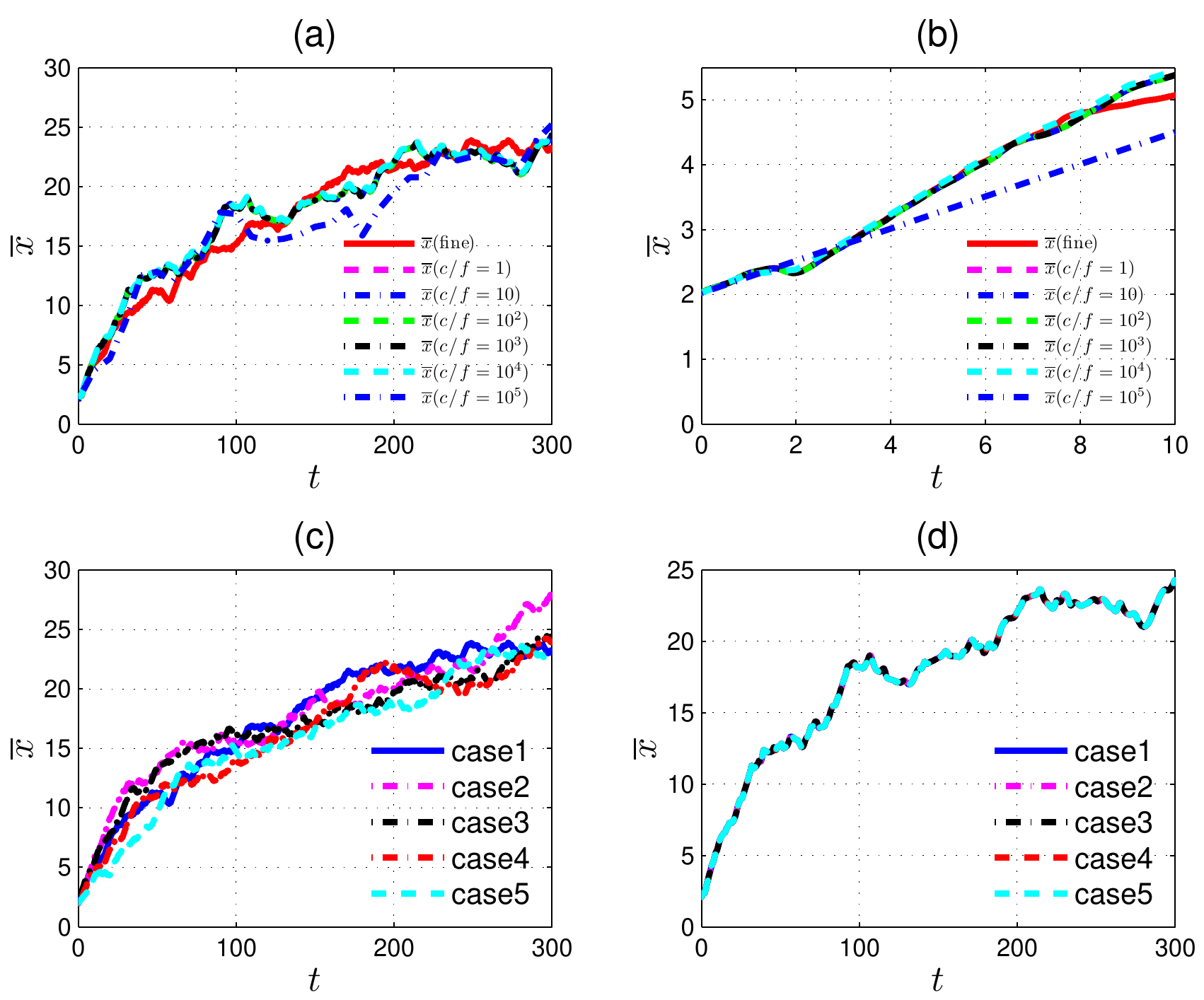}
\makeatletter
\def\@captype{figure}
\makeatother
\caption{(a) Comparison of coarse and exact response for different coarse-to-fine time-step ratios; (b) Magnified version of (a) from $t$=0 to 10; (c) $\overline{x}(t)$ from fine theory with different fine ICs, where the $l^2$-norm of any two corresponding coarse ICs is within $10^{-3}$; (d) $\overline{x}(t)$ from coarse theory, where the coarse ICs are the same as (c)}\label{Lorenz_fig3}
\end{center}

\begin{center}
\includegraphics [width=.4\columnwidth]{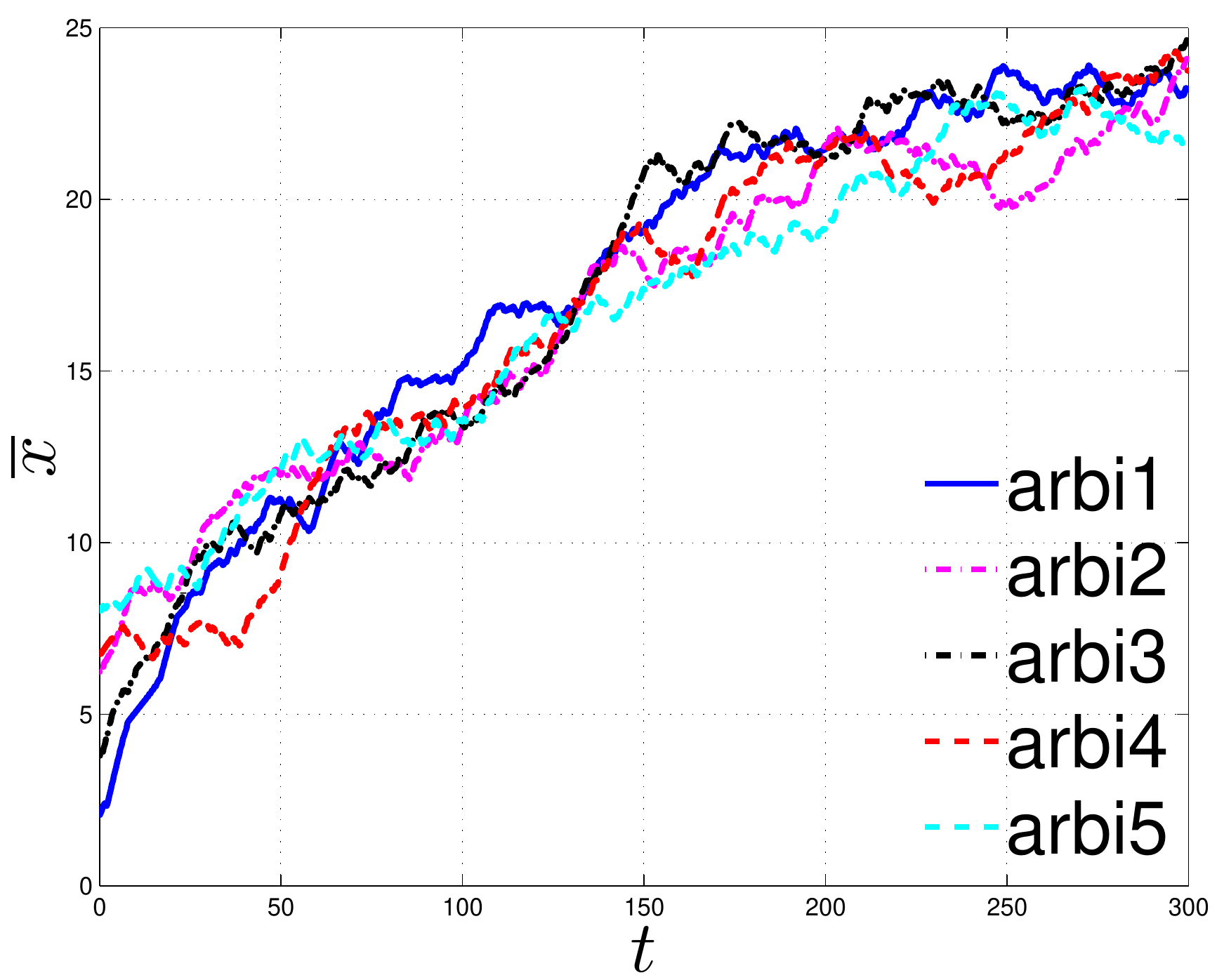}
\makeatletter
\def\@captype{figure}
\makeatother
\caption{$\overline{x}(t)$ from fine theory with arbitrary fine ICs}\label{Lorenz_fig4}
\end{center}

\noindent
In conclusion, we succeed in obtaining consistent coarse response from solving the coarse theory with much larger time-step. We learn that by taking a large enough averaging period, we can obtain a coarse dynamics that is insensitive to fine ICs, while retrieving some key features at the coarse level; on the other hand, the coarse dynamics is not aimed to capture the chaotic feature and high frequency dynamics at the fine scale. With small differences of coarse ICs, the coarse responses show no difference. Both of these properties are desirable since the dynamics at much larger time scale than that of fine dynamics is of interest to us.\\

\noindent
{\bfseries 3. Approximation of High Frequency Oscillatory Response}\\

\noindent
The problems of interest here are singularly perturbed systems that display highly oscillatory response and whose limiting behavior often retains fast oscillatory response. Many physical problems can be expected to fall in this category. For instance, in MD simulations where external loading is applied to the system, the loading rate is usually at least $10^9$ times slower than the fundamental frequency of atomic vibrations. This large separation of time-scales leads to singular perturbation in the fine system~(\ref{eq:org}) since $1/\mbox{T}$ becomes a small number. In this study, we are interested in understanding the evolution of \emph{time averages} of phase functions of the fast flow on the slow time scale of loading. Two approaches are developed to approximate the time averaged dynamics. The first approach utilizes the Tikhonov method~\cite{tikhonov1952systems} of differential-algebraic equations (DAE). The second approach is a natural, practical extension of the method of averaging proposed in~\cite{artstein1996singularly, artstein2007slow}. We construct the evolutionary equations of time averaged quantities following the work of~\cite{slemrod2010time} and refer to this collection of theoretical and numerical ideas as `Practical Time Averaging'.  \\

\noindent
3.1 THE TIKHONOV APPROACH \\

\noindent
The system we are concerned with has the same form
\begin{equation}
\begin{aligned}
&\frac{df}{dt}=H(f,I) \\
&\frac{dI}{dt}=\frac{1}{\mbox{T}}L(I) \\
&c(t)=\frac{1}{\tau}\int_t^{t+\tau}\Lambda(f(p))dp,
\end{aligned}
\end{equation}
where $f\in{\Bbb{R}^N}$, $I\in{\Bbb{R}^n}$ and $c\in{\Bbb{R}^m}$. Recall that if $\omega$ is the loading frequency, then $\mbox{T}\rightarrow\infty$ as $\omega \rightarrow 0$. Define a new time scale $s:={t}/{\mbox{T}}$ and denote $\varepsilon:={1}/{\mbox{T}}$. Note that since T is dimensionless, $s$ has the unit of time and $\varepsilon$ is a dimensionless number. For the functions $f$ and $I$ of the fast time variable $t$, we define the functions $\tilde{f}$ and $\tilde{I}$ of the slow time variable $s$ as follows:
\begin{equation}
\tilde{f}(s)=f(s\mbox{T}),  \;\;\;\; \tilde{I}(s)=I(s\mbox{T}).
\end{equation}
Then, on the slow time scale the system reads
\begin{equation}
\begin{aligned}
&\varepsilon\frac{d\tilde{f}}{ds}=H(\tilde{f},\tilde{L})  \\
&\frac{d\tilde{I}}{ds}=L(\tilde{I}).
\end{aligned}
\end{equation}
The coarse variable is defined as
\begin{equation}
\tilde{c}(s)=c(s\mbox{T})=\frac{1}{\tau}\int_{s\mbox{T}}^{s\mbox{T}+\tau}\Lambda(f(p))dp,
\end{equation}
and we assume that $\tau=\kappa/\omega$ with $\kappa$ a non-dimensional constant and $0<\kappa<1$. Introducing the change of variables $r=p/\mbox{T}$, we have
\begin{equation}
\tilde{c}(s)=\frac{1}{\lambda}\int_s^{s+\lambda}\Lambda(\tilde{f}(r))dr,
\end{equation}
where $\lambda=\tau/\mbox{T}$. Therefore, $\lambda=\kappa/(\omega\mbox{T})=\kappa/\omega_f$ is a constant as $\varepsilon\rightarrow0$, with physical dimensions of time.\par
Removing \emph{all overhead tildes}, we have the slow-fast system, posed on the slow time scale:
\begin{equation}\label{slow_fast}
\begin{aligned}
&\varepsilon\frac{df}{ds}=H(f,I) \\
&\frac{dI}{ds}=L(I) \\
&c(s)=\frac{1}{\lambda}\int_{s}^{s+\lambda}\Lambda(f(r))dr.
\end{aligned}
\end{equation}
Let $\hat{f}$ be the solution of the following differential-algebraic system obtained from~(\ref{slow_fast}) when $\varepsilon$ is set to be equal to 0, namely
\begin{equation}\label{DAE}
\begin{aligned}
&0=H(\hat{f},I) \\
&\frac{dI}{ds}=L(I).
\end{aligned}
\end{equation}
Then by using the chain rule, we obtain
\begin{equation}
\frac{\partial H}{\partial \hat{f}}\frac{d\hat{f}}{ds}=-\frac{\partial H}{\partial I}\frac{dI}{ds}.
\end{equation}
So the new evolution equation (rate-independent) for the fine variables is
\begin{equation}
\frac{d\hat{f}}{ds}=-\left(\frac{\partial H}{\partial \hat{f}}\right)^{-1}\frac{\partial H}{\partial I}\frac{dI}{ds}.
\end{equation}
The new system becomes
\begin{equation}
\begin{aligned}
&\frac{d\hat{f}}{ds}=-\left(\frac{\partial H}{\partial\hat{f}}\right)^{-1}\frac{\partial H}{\partial I}L \\
&\frac{dI}{ds}=L(I) \\
&c(s)=\frac{1}{\lambda}\int_{s}^{s+\lambda}\Lambda(\hat{f}(r))dr
\end{aligned}\label{Approx_dynamics}
\end{equation}
with initial conditions on $\hat{f}$ satisfying $H(\hat{f},I)=0$. \par
A classical result of Tikhonov~\cite{tikhonov1952systems} states the solutions of~(\ref{slow_fast}) converge, as $\varepsilon$ tends to zero, to the solution of the differential-algebraic system~(\ref{DAE}) under appropriate assumptions. A crucial condition for this result to hold is that for each fixed $I$, solutions of the differential equation ${df}/{dt}=H(f,I)$ should converge, as $t\rightarrow\infty$, to a stationary point $\hat{f}(I)$ such that $H(\hat{f}(I),I)=0$. However, a more general theory given by~\cite{artstein1996singularly, artstein2002singularly} shows that this condition may not be always satisfied. Instead, the convergence to the limit is in a weaker sense (in the sense of Young measures), namely, one gets information only about the probability measure of the fast flow. Our work is to test whether the slow time evolution of the first moment of the probability measure of the fast flow parametrized by $I(s)$ corresponds in any sense to the solution of the DAE of the Tikhonov approach. Realizing that~(\ref{DAE}) and~(\ref{Approx_dynamics}) are equivalent given the ODE starts from an equilibrium point of the fast dynamics, we refer to the ODE system~(\ref{Approx_dynamics}) as the DAE dynamics of the slow-fast system~(\ref{slow_fast}). The system, in general, is not well-defined when ${\partial H}/{\partial \hat{f}}$ is singular. However, well-posed evolution can exist in many situations. At the practical level, this may be solved by trying different numerical algorithms (for example, using implicit Euler backward method on the points near the singularity), but we avoid an extensive consideration of this fact in the following model problems. \\

\noindent
3.1.1 Singularly Perturbed System with Measure-valued Limit \\

\noindent
We discuss a specific realization of Example 5.1 in~\cite{artstein2002singularly}. The slow-fast system is given as
\begin{equation}\label{arstein52}
\begin{aligned}
&\varepsilon\frac{y_1}{ds}=y_2 \\
&\varepsilon\frac{y_2}{ds}=g(x,y_1,y_2) \\
&\frac{dx}{ds}=y_1,
\end{aligned}
\end{equation}
with state variables $y_1$, $y_2$ and $x$. The function $g(x,y_1,y_2)$ is defined as
\begin{equation}
g(x,y_1,y_2)=\alpha(x)(1-A(x)(y_1-z_{up}(x))^2)y_2-(y_1-z_{up}(x)),
\end{equation}
where we set $\alpha(x)=(x-3^{-\frac{3}{2}}2)^2$, $A(x)=1/(x-3^{-\frac{3}{2}}2)$ and
\begin{equation}
z_{up}(x)=\{\hat{y}_1:-x+\hat{y}_1-\hat{y}_1^3=0\; \forall\;x\in{\Bbb{R}}\;\; \mbox{and}\;\; \hat{y}_1>3^{-\frac{1}{2}}\}.
\end{equation}
for the dynamics at the upper branch (see Figure~\ref{Approx_Dynamics_Arstein}). For the lower branch, we substitute $\alpha(x)=(x-(-3^{-\frac{3}{2}}2))^2$, $A(x)=1/(x-(-3^{-\frac{3}{2}}2))$ and replace $z_{up}(x)$ with $z_{lo}(x)$ which is defined as
\begin{equation}
z_{lo}(x)=\{\hat{y}_1:-x+\hat{y}_1-\hat{y}_1^3=0\; \forall\;x\in{\Bbb{R}}\;\; \mbox{and}\;\; \hat{y}_1<-3^{-\frac{1}{2}}\}.
\end{equation}
Thus the corresponding differential-algebraic system for~(\ref{arstein52}) has the form
\begin{equation}\label{arstein_DAE}
\begin{aligned}
&0=\hat{y}_2 \\
&0=-x+\hat{y}_1-\hat{y}_1^3 \\
&\frac{dx}{ds}=\hat{y}_1.
\end{aligned}
\end{equation}
However, for $\varepsilon\rightarrow0$ the limit behavior of~(\ref{arstein52}) does not follow the DAE~(\ref{arstein_DAE}). Instead, for each fixed $x$, the fast motion will converge to periodic orbits (also referred to as limit cycles) around either the upper or lower branch of the equilibria profile given by~(\ref{arstein_DAE}), as shown in Figure~\ref{Approx_Dynamics_Arstein}. The blue curve exhibits the trajectory from running the slow-fast system~(\ref{arstein52}), plotted on $x$-$y_1$ phase space. Meanwhile, there are fast motions around $y_2=0$ (not shown) during evolution. The equilibria profile is shown in red with $\hat{y}_2$ equal to zero constantly. We start the dynamics at $x=-1.875$. The initial values of $(y_1,y_2)$ are set to be (1.5,0). However, in fact, starting from any values of $(y_1,y_2)$ with $y_1>z_{up}(x)$ will lead the fast motions to converge to the limit cycle around the equilibrium point within a short time interval with $x$ hardly changing. The trajectory continues its fast motion in $y_1$-$y_2$ plane close to the limit cycles (even in the limit $\varepsilon\rightarrow0$), never reaching the slow manifold defined by $g(x,y_1,y_2)= 0$ with a slow movement along the positive $x$ direction till $x=3^{\frac{3}{2}}2$. Similarly, when starting from a point on the lower branch, say $x=1.875$, the fast solution will converge to the limit cycle shortly following a slow drift along the negative $x$ direction till the point $x=-3^{\frac{3}{2}}2$. \par
The DAE dynamics for~(\ref{arstein52}) according to~(\ref{Approx_dynamics}) is
\begin{equation}
\begin{aligned}
&\frac{d\hat{y}_1}{ds}=-\left(\frac{\partial g}{\partial x}\right)\bigg/
\left(\frac{\partial g}{\partial \hat{y}_1}\right)\hat{y}_1   \\
&\frac{d\hat{y}_2}{ds}=0 \\
&\frac{dx}{ds}=\hat{y}_1.
\end{aligned}
\end{equation}
The numerical result is shown in Figure~\ref{Approx_Dynamics_Arstein52}. Initial conditions are chosen to satisfy $H=0$. There are two singularity points: $(x,\hat{y}_1,\hat{y}_2)=(3^{\frac{3}{2}}2, 3^{-\frac{1}{2}}, 0)$ and $(x,\hat{y}_1,\hat{y}_2)=(-3^{\frac{3}{2}}2, -3^{-\frac{1}{2}}, 0)$. We start the simulation from the point $(-1.875, 1.5, 0)$ on the upper branch until the trajectory reaches the singularity point $(3^{\frac{3}{2}}2, 3^{-\frac{1}{2}}, 0)$, and then repeat the same procedure on the lower branch till the point $(-3^{\frac{3}{2}}2, -3^{-\frac{1}{2}}, 0)$. $\hat{y}_2$ is constant zero in both cases. It can be seen that the time averaged behavior matches well with the DAE dynamics. The two trajectories are almost overlapping. Although the DAE is not the appropriate limit solutions of the system, it succeeds in capturing the time averaged behavior of the limit dynamics.\\
\begin{center}
\includegraphics [width=.4\columnwidth]{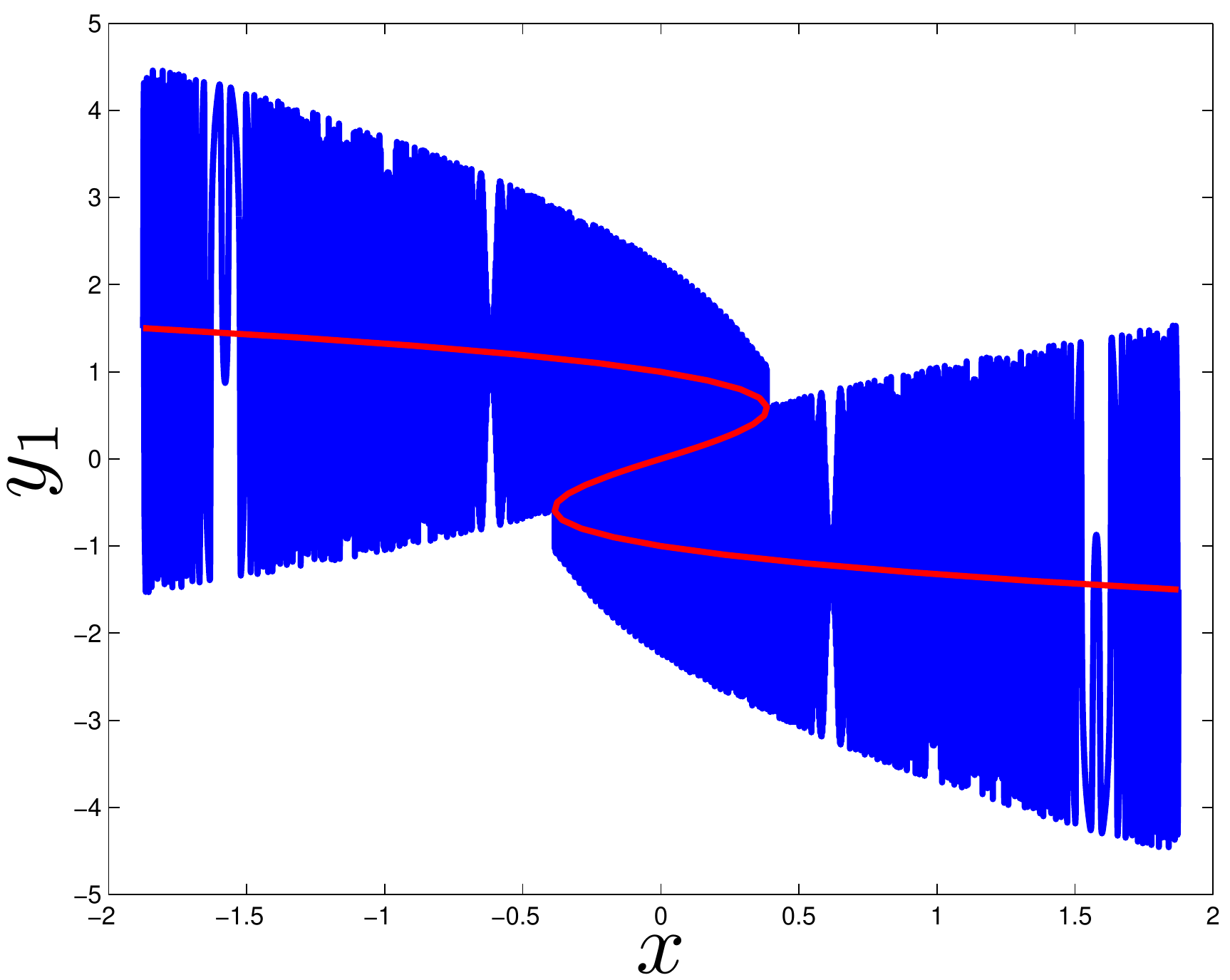}
\makeatletter
\def\@captype{figure}
\makeatother
\caption{\small Diagram of equilibria and limit cycles. Red curve indicates the equilibria of the fast system and blue curve shows the limit cycles generated from the slow-fast system in the limit $\varepsilon \rightarrow 0$.} \label{Approx_Dynamics_Arstein}
\end{center}

\begin{center}
\includegraphics [width=.4\columnwidth]{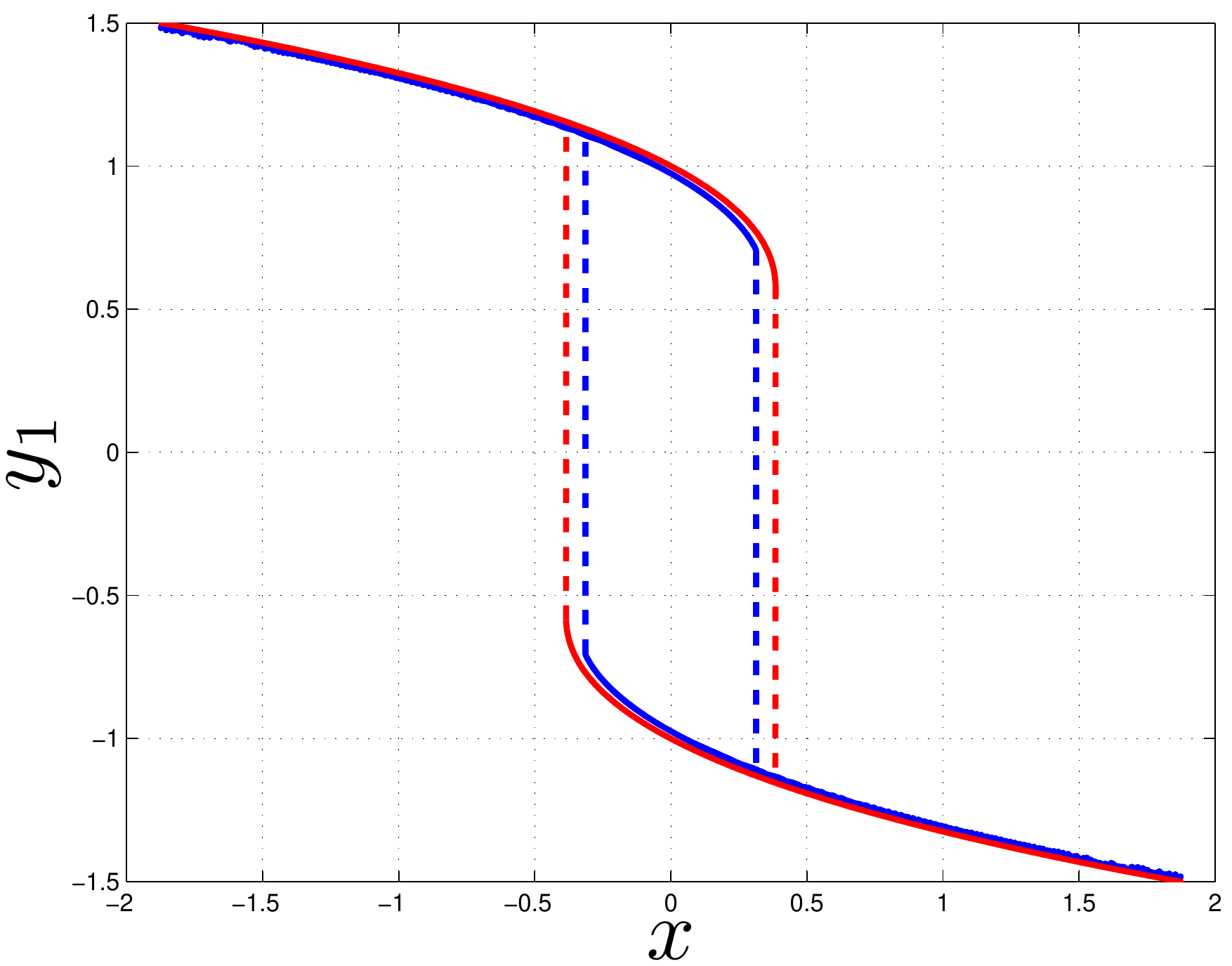}
\makeatletter
\def\@captype{figure}
\makeatother
\caption{\small Comparison of time averaged dynamics and DAE dynamics. Blue curve indicates $\overline{y}_1$ vs. $x$; red curve indicates $\hat{y}_1$ vs. $x$. Solid line indicates stable dynamics; dash line indicates fast motion of limit dynamics.} \label{Approx_Dynamics_Arstein52}
\end{center}

\noindent
3.1.2 `Oscillatory Forced' Lorenz System \\

\noindent
The original Lorenz system is a three-dimensional chaotic system. The trajectories tend to converge on a topologically complicated, but bounded set that can be enveloped within another bounded and contiguous region of 3D phase space. We apply the same strategy as in Section 2.2 to make the averaged $x$ variable evolve sinusoidally with loading period while keeping the high frequency oscillation. By doing so, we assure that the trajectories of the `Forced' Lorenz system will remain bounded (under the condition that the loading is bounded) and exhibit chaotic patterns along a slow sinusoid. Meanwhile, the `Forced' Lorenz system will show an explicit separation of slow and fast motions. Towards this goal, we define
\begin{equation}\label{defx2}
\begin{aligned}
&\tilde{x}(t)=x(t)+L(t) \\
&\mbox{where} \;\; \dot{L}=\omega{g}\;\; \mbox{and} \;\; \dot{g}=-\omega{L}.
\end{aligned}
\end{equation}
Obviously $\omega$ indicates the loading frequency. From~(\ref{defx2}), we have
\begin{equation}
x(t)=\tilde{x}(t)-L(t)
\end{equation}
and
\begin{equation}\label{defx3}
\dot{\tilde{x}}=\dot{x}+\dot{L}.
\end{equation}
Substitute~(\ref{defx2}) and~(\ref{defx3}) into the original Lorenz system~(\ref{eq:lorenz}), and remove the overhead tilde notation,
\begin{equation}
\begin{aligned}
&\dot{x}=\sigma(y-x)+\sigma{L}+\omega{g} \\
&\dot{y}=\gamma{x}-y-xz-\gamma{L}+Lz \\
&\dot{z}=xy-\beta{z}-Ly \\
&\dot{L}=\omega{g} \\
&\dot{g}=-\omega{L}.
\end{aligned}
\end{equation}
Assume $s=\omega{t}$, the system on the slow time scale takes the form
\begin{equation}\label{org_Lorenz}
\begin{aligned}
&\omega\frac{dx}{ds}=\sigma(y-x)+\sigma{L}+\omega{g}  \\
&\omega\frac{dy}{ds}=\gamma{x}-y-xz-\gamma{L}+Lz  \\
&\omega\frac{dz}{ds}=xy-\beta{z}-Ly  \\
&\frac{dL}{ds}=g  \\
&\frac{dg}{ds}=-L.
\end{aligned}
\end{equation}
For $\omega\rightarrow{0}$, the equilibrium points of the fast system at $s=0$ are denoted as $\mbox{EP1}=(8,8,24)$, $\mbox{EP2}=(-8,-8,24)$ and $\mbox{EP3}=(0,0,0)$. \par
The DAE dynamics following~(\ref{DAE})-(\ref{Approx_dynamics}) is
\begin{equation}\label{DAE1}
\begin{aligned}
&\frac{d\hat{x}}{ds}=\frac{\omega{L}[(L-\hat{x})(\hat{x}-L)-\beta]}{\sigma(-2L\hat{x}+L^2+\hat{x}^2-\beta\gamma+\beta\hat{z}+\beta-L\hat{y}+\hat{x}\hat{y})}+g \\
&\frac{d\hat{y}}{ds}=-\frac{\omega{L}(L\hat{y}-\beta\hat{z}+\beta\gamma-\hat{x}\hat{y})}{\sigma(-2L\hat{x}+L^2+\hat{x}^2-\beta\gamma+\beta\hat{z}+\beta-L\hat{y}+\hat{x}\hat{y})} \\
&\frac{d\hat{z}}{ds}=-\frac{\omega{L}(\hat{y}+\gamma\hat{x}-\gamma{L}-\hat{z}\hat{x}+\hat{z}L)}{\sigma(-2L\hat{x}+L^2+\hat{x}^2-\beta\gamma+\beta\hat{z}+\beta-L\hat{y}+\hat{x}\hat{y})} \\
&\frac{dL}{ds}=g \\
&\frac{dg}{ds}=-L.
\end{aligned}
\end{equation}
Since $\omega$ exists in the right hand side terms in~(\ref{org_Lorenz}), it remains in the DAE dynamics. With further assumption that $\omega\rightarrow{0}$ in~(\ref{DAE1}), we obtain a simple form of the DAE dynamics:
\begin{equation}\label{Lorenz_approx}
\begin{aligned}
&\frac{d\hat{x}}{ds}=g \\
&\frac{d\hat{y}}{ds}=0 \\
&\frac{d\hat{z}}{ds}=0 \\
&\frac{dL}{ds}=g \\
&\frac{dg}{ds}=-L.
\end{aligned}
\end{equation}
Figure~\ref{Approx_Dynamics_Lorenz} shows the fine variables and their time averages from evolving the system~(\ref{org_Lorenz}). The evolution of $x$ is a superposition of the high frequency oscillation with the periodic loading whose time average coincides with the sinusoidal loading. The time averages of $y$ and $z$ are nearly constants, despite the oscillatory motion on the fast time scale. Figure~\ref{Approx_Dynamics_Lorenz_bar}(a),~\ref{Approx_Dynamics_Lorenz_bar}(b) and~\ref{Approx_Dynamics_Lorenz_bar}(c) display the time averages of $x$, $y$ and $z$, respectively, starting from different initial conditions. It is noted that the averaged behavior, or the dynamics on the slow time scale, is insensitive to the fine initial conditions (for the cases tested). The trajectories that start from different points converge to a \emph{unique} time averaged dynamics (this statement is limited to the numerical tests we did) and follow the time averaged dynamics consistently.\par
However, the DAE dynamics~(\ref{Lorenz_approx}) does not give consistent results. The trajectories starting from equilibria do not coincide with the time averaged dynamics, as seen in Figure~\ref{Approx_Dynamics_Lorenz_hat}. Figure~\ref{Approx_Dynamics_Lorenz_hat}(a), (b) and (c) present the results from fine dynamics and DAE dynamics of the variables $x$, $y$ and $z$, respectively. The black curve denotes the time average from the fine theory. Although the trajectories evolve with the same pattern, they do not agree with the time averaged behavior since the initial conditions corresponding to equilibria of fast dynamics do not coincide with the averaged dynamics. The last test is to set initial conditions for the DAE dynamics as the average of $(\overline{x},\overline{y},\overline{z})$ from previously assigned fine initial conditions. In this case, $(\hat{x},\hat{y},\hat{z})$ are consistent with $(\overline{x},\overline{y},\overline{z})$ (see the magenta curve in Figure~\ref{Approx_Dynamics_Lorenz_hat}).\par
There is an analogue between the `Forced' Lorenz dynamics and Example 5.1 discussed in Section 3.1.1 as can be seen from a plot like Figure~\ref{Approx_Dynamics_Lorenz}(a). The fast motion occurs in $(x,y,z)$ space (the dynamics in the $y$ and $z$ directions is not shown) with $x$ oscillating around the slow sinusoid and $(y,z)$ moving around a fixed point close to (0, 20.6).\\

\begin{center}
\includegraphics [width=.8\columnwidth]{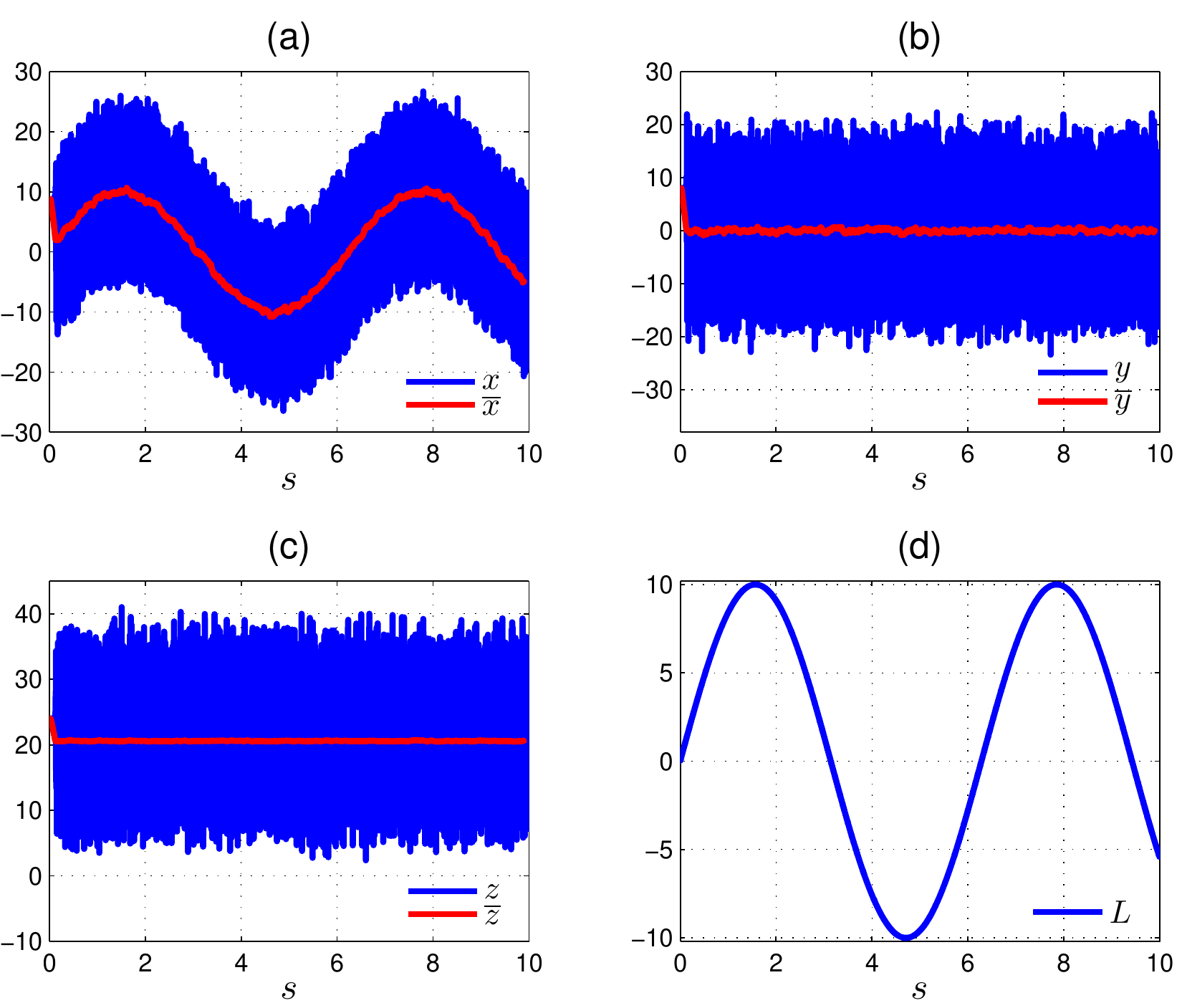}
\makeatletter
\def\@captype{figure}
\makeatother
\caption{\small Comparison of fine evolution and corresponding time averaged dynamics} \label{Approx_Dynamics_Lorenz}
\end{center}

\begin{center}
\includegraphics [width=.8\columnwidth]{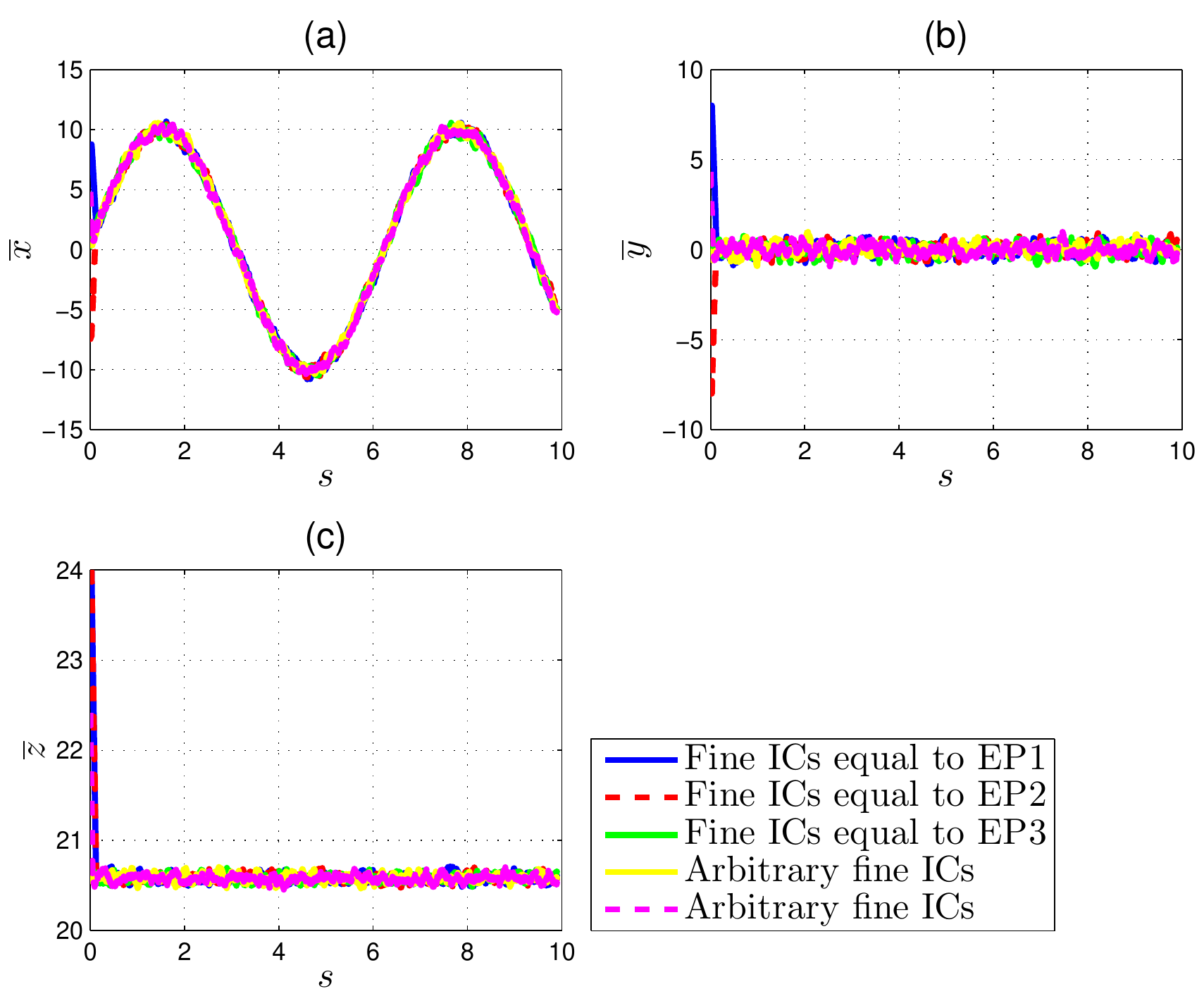}
\makeatletter
\def\@captype{figure}
\makeatother
\caption{\small Time averaged dynamics starting from various initial conditions} \label{Approx_Dynamics_Lorenz_bar}
\end{center}

\begin{center}
\includegraphics [width=.8\columnwidth]{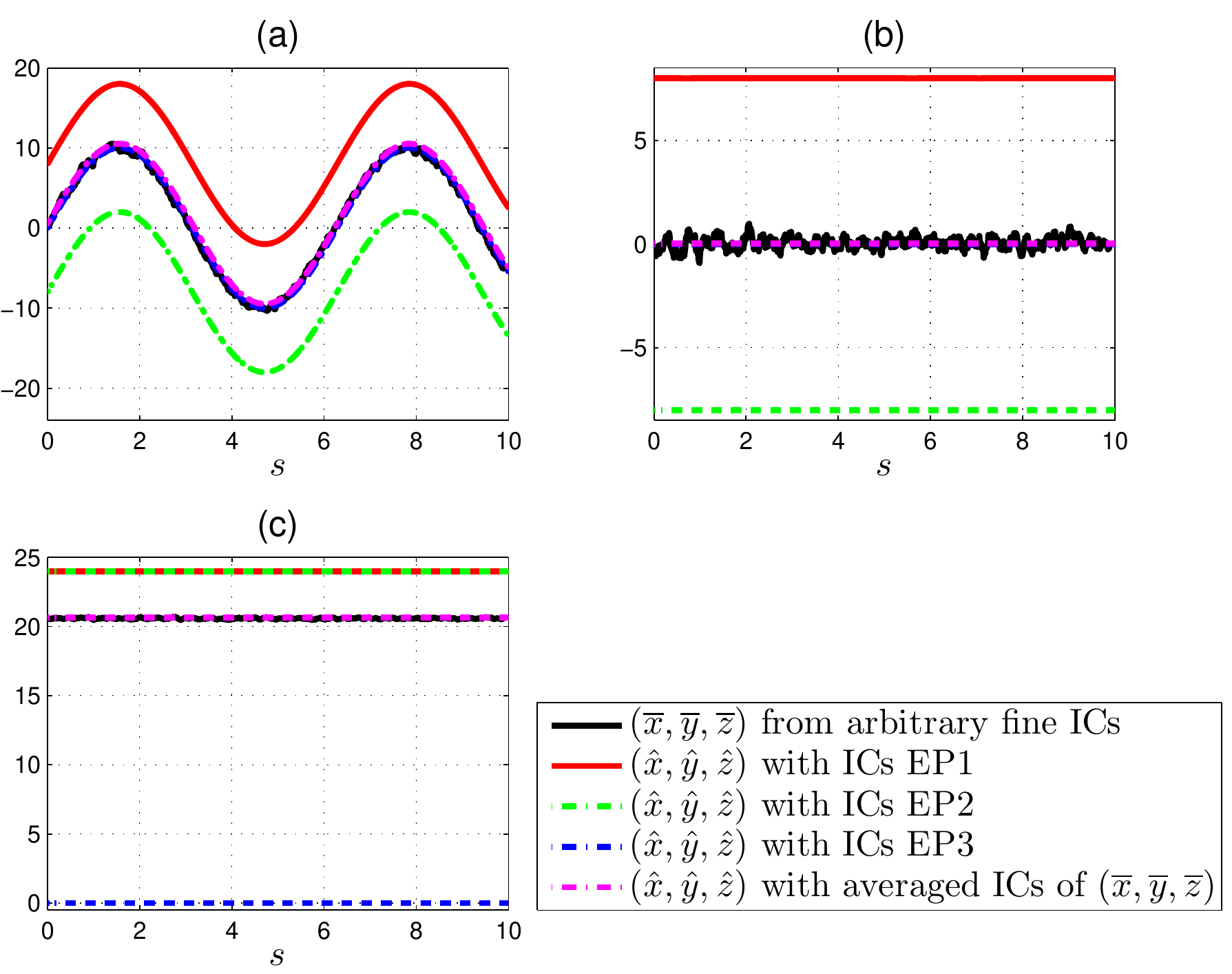}
\makeatletter
\def\@captype{figure}
\makeatother
\caption{\small Comparison of time averaged dynamics from fine and DAE dynamics starting from various initial conditions}\label{Approx_Dynamics_Lorenz_hat}
\end{center}

\clearpage
\noindent
3.2 PRACTICAL TIME AVERAGING \\

\noindent
Different from the Tikhonov approach, but including it when appropriate, the work of~\cite{artstein2007slow} building on earlier work of~\cite{artstein1996singularly} allows dealing with limits of fast flows that are not necessarily equilibria. The theory accounts for the limit behavior of fast flows that may be supported on $\omega$-limit sets of complicated topology. We state the relevant part of the theorem below, following the terminology presented in~\cite{slemrod2010time}:\par
Suppose the dynamical system of interest, written on the slow time-scale, is~(\ref{slow_fast}). Let $f^{\varepsilon,I,f_0}$ represent solutions of the fast system
\begin{equation}\label{PTA_fast}
\varepsilon\frac{df^{\varepsilon,I,f_0}}{ds}=H(f^{\varepsilon,I,f_0},I)\quad ; \quad f^{\varepsilon,I,f_0}(0)=f_0
\end{equation}
with $I$ fixed. Rewrite the fast system with the introduction of the fast time scale $t=s/\varepsilon$,
\begin{equation}\label{fast_system}
\frac{df^{(0),I,f_0}}{dt}(t)=H(f^{(0),I,f_0}(t),I)\; ;\;f^{(0),I,f_0}(0)=f_0
\end{equation}
where
\begin{equation}\label{change_timescales}
f^{(0),I,f_0}(t)=f^{(0),I,f_0}\left(\frac{s}{\varepsilon}\right)=f^{\varepsilon,I,f_0}(s).
\end{equation}
Assume that for $0\leq{t}<\infty$, $f^{\varepsilon,I,f_0}$ lie in a bounded set for each $I$. Hence, on $[0,\mathsf{T}]$ we will have for any continuous function $F$
\begin{equation}\label{Young_measure}
F(f^{\varepsilon,I,f_0}(s))\rightarrow\int_{\Bbb{R}^N}F(\gamma)\mu_{s,I,f_0}(\gamma)d\gamma\;\; \mbox{weak}\;\star\;L^{\infty}[0,\mathsf{T}],
\end{equation}
where $\mu_{s,I,f_0}$ is the Young measure generated by the sequence $f^{\varepsilon,I,f_0}\in L^{\infty}(0,\mathsf{T};\Bbb{R}^N)$. Moreover, the theory proves that $\mu_{s,I,f_0}$ is an invariant measure supported on the $\omega$-limit set of $f^{(0),I,f_0}$. \par
Note that $f^{\varepsilon,I,f_0}$ is obtained from running the fast system~(\ref{fast_system}) for a sufficiently long period of time, starting from $f_0$ while keeping $I$ fixed; The function $F$ in~(\ref{Young_measure}) can be linear or nonlinear. In situation when the $\omega$-limit set of the fast dynamics is not an equilibrium point, the Tikhonov approach may not in general be expected to approximate the time evolution of the first-moment of nonlinear functions of state. Furthermore, the theory proves that coarse evolution equations can be written down for state variables that are `slow', i.e. orthogonal to the fast flow in a precise sense, but does not provide a prescription for generating such `orthogonal observables'. \par
We follow the work of~\cite{slemrod2010time}, which provides a straightforward way of generating a large class of slow variables via time averaging that are shown to be `orthogonal observables' (for a corresponding infinite dimensional dynamics) for which a slow dynamics can be written down. We refer to this approach and the numerical approximation ideas we develop in this paper as `Practical Time Averaging'. \par
Recall the coarse variable defined on the slow time scale~(\ref{slow_fast}) and write down its evolution:
\begin{equation}\label{coarse_var}
c(s)=\frac{1}{\lambda}\int_s^{s+\lambda}\Lambda(f(r))dr\Longrightarrow\frac{dc}{ds}(s)=\frac{1}{\lambda}\left[\Lambda(f(s+\lambda))-\Lambda(f(s))\right].
\end{equation}
Let us define the sequence (on $\varepsilon$) of smooth functions $c^{\varepsilon,f_0}$ such that
\begin{equation}\label{avg_inv}
c^{\varepsilon,f_0}(s)=\frac{1}{\lambda}\int_s^{s+\lambda}\Lambda(f^{\varepsilon,I(r),f_0}(r))d{r},
\end{equation}
where we follow the notation in~(\ref{PTA_fast}) but now $I$ is a function of the slow time.
Then the evolution equation is
\begin{equation}\label{coarse_invariant}
\frac{d{c^{\varepsilon,f_0}}}{ds}(s)=\frac{1}{\lambda}\left[\Lambda(f^{\varepsilon,I(s+\lambda),f_0}(s+\lambda))-\Lambda(f^{\varepsilon,I(s),f_0}(s))\right].
\end{equation}
Let $\overline{c}^{f_0}$ denote the weak limit of the sequence $c^{\varepsilon,f_0}$, and using~(\ref{Young_measure}), we have
\begin{equation}\label{limit_invariant}
\frac{d\overline{c}^{f_0}}{ds}=\frac{1}{\lambda}\left(\int_{\Bbb{R}^N}\Lambda(\gamma)\mu_{s+\lambda,I(s+\lambda),f_0}(\gamma)d\gamma-\int_{\Bbb{R}^N}\Lambda(\gamma)\mu_{s,I(s),f_0}(\gamma)d\gamma \right)
\end{equation}
in the weak (variational) sense. (\ref{limit_invariant}) can be considered as the instantaneous evolution equation of the limit slow variable $\overline{c}^{f_0}$ of the sequence $c^{\varepsilon, f_0}$. It should be noted that if any arbitrary instantaneous function, say $A(f)$, is taken as the coarse variable, we have
\begin{equation}\label{instantaneous}
\frac{dA(f)}{ds}(s)=\frac{1}{\varepsilon}{\nabla}_f{A}\cdot{H}.
\end{equation}
Since ${1}/{\varepsilon}$ is built in the evolution equation, the Young measure limit is usually not applicable unless the function $A(f)$ satisfies $\nabla_fA\cdot{H}=0$. However, our definition of coarse variables allows constructing coarse evolution equation for more general kinds of averaged phase functions. For example, to approximate the change of temperature of a MD system, where temperature can be interpreted as a function of particle velocities but normally does not satisfy $\nabla_fA\cdot{H}=0$, we can define the time average of temperature and follow the procedure from~(\ref{coarse_var}) to~(\ref{limit_invariant}) to write down the coarse evolution equation. \par
To numerically approximate the coarse dynamics~(\ref{limit_invariant}), we can apply a three-stage procedure suggested in~\cite{artstein2009analysis}:
\begin{enumerate}
\item Evolve the fast system~(\ref{fast_system}) with slow variables fixed at time step $s$ and $s+\lambda$, respectively to generate the Young measures as averages of $M$ Dirac masses at $M$ values of $f$, i.e.
\begin{equation}
\begin{aligned}
&\frac{1}{\lambda}\int_{\Bbb{R}^N}\Lambda(\gamma)\mu_{s,I(s),f_0}(\gamma)d\gamma\approx\frac{1}{M}\sum_{i=1}^{M}\Lambda(f^{(0),I(s),f_0})  \\
&\frac{1}{\lambda}\int_{\Bbb{R}^N}\Lambda(\gamma)\mu_{s+\lambda,I(s+\lambda),f_0}(\gamma)d\gamma\approx\frac{1}{M}\sum_{i=1}^{M}\Lambda(f^{(0),I(s+\lambda),f_0}),
\end{aligned}
\end{equation}
where $M$ is chosen to be large enough for the averages to converge.
\item Extrapolate to obtain the coarse variable at the next coarse step $s+T$ using the Forward Euler method,
\begin{equation}\label{step3}
\overline{c}^{f_0}(s+T)\approx\overline{c}^{f_0}(s)+\frac{T}{\lambda}\left[\int_{\Bbb{R}^N}\Lambda(\gamma)\mu_{s+\lambda,I(s+\lambda),f_0}(\gamma)d\gamma-\int_{\Bbb{R}^N}\Lambda(\gamma)\mu_{s,I(s),f_0}(\gamma)d\gamma\right],
\end{equation}
where $T$ is the coarse time step. Usually $T$ is set to be much larger than $\lambda$. $I(s+T)$ can be obtained from evolving~$(\ref{slow_fast})_2$ solely.
\item Reconstruct to find fine data that matches $\overline{c}^{f_0}(s+T)$ and $I(s+T)$.
\item Repeat.
\end{enumerate}\par
The primary shortcoming of all of this theory as acknowledged in~\cite{artstein1996singularly, artstein2007slow} and carries over for the model with time averages in~\cite{slemrod2010time} is that the theory is silent about fine initial conditions required to generate the Young measures at the discrete slow time-levels. This is a serious practical issue and renders $f_0$ as an essentially free parameter whose choice, nevertheless, affects results and has to be made on an ad-hoc basis based on numerical experimentation. We mention here that the unsubstantiated assumption of ergodicity is often made to sweep this important issue under the rug.\par
Given the above uncertainty and the need to use it twice to follow the evolution prescribed by~(\ref{step3}), an alternative procedure is the following. Based on~(\ref{change_timescales}), the sequence of functions $c^{\varepsilon,f_0}(s)$ can be written as
\begin{equation}\label{change_lhs_timescales}
c^{\varepsilon,f_0}(s)=c^{(0),f_0}\left(\frac{s}{\varepsilon}\right)=c^{(0),f_0}(t).
\end{equation}
Using the substitutions $p={r}/{\varepsilon}$, $\tau=\lambda/\varepsilon$ and~(\ref{change_timescales}), the right-hand-side term in~(\ref{avg_inv}) can be interpreted on the fast time scale as
\begin{equation}\label{change_rhs_timescales}
\frac{1}{\lambda}\int_s^{s+\lambda}\Lambda(f^{\varepsilon,I(r),f_0}(r))d{r}\approx\frac{1}{\tau}\int_t^{t+\tau}\Lambda(f^{(0),I(t),f_0}(p))d{p},
\end{equation}
where we make the realistic assumption that the loads $I$ do not vary appreciably in the (fast-time) interval $[t, t+\tau]$. From~(\ref{change_lhs_timescales}) and~(\ref{change_rhs_timescales}), we have
\begin{equation}\label{alternative_coarse}
c^{(0),f_0}(t)=\frac{1}{\tau}\int_t^{t+\tau}\Lambda(f^{(0),I(t),f_0}(p))d{p},
\end{equation}
which is essentially~(\ref{eq:def_cos}), of course.
We are interested in $c^{\varepsilon,f_0}$ as $\varepsilon \rightarrow 0$, which by definition, is the time average of functions of the fast variables as $\tau \rightarrow \infty$ when the slow variable $I$ is fixed. We simply follow~(\ref{alternative_coarse}) to calculate this average at each coarse time step $T$. Theoretically, we need to run the fast system for an infinite period of time. However, for the numerical tests conducted here we have seen that the average converges quickly, thus taking a reasonable period of time is sufficient.\par
\begin{figure}
\centering
\includegraphics [width=.4\columnwidth]{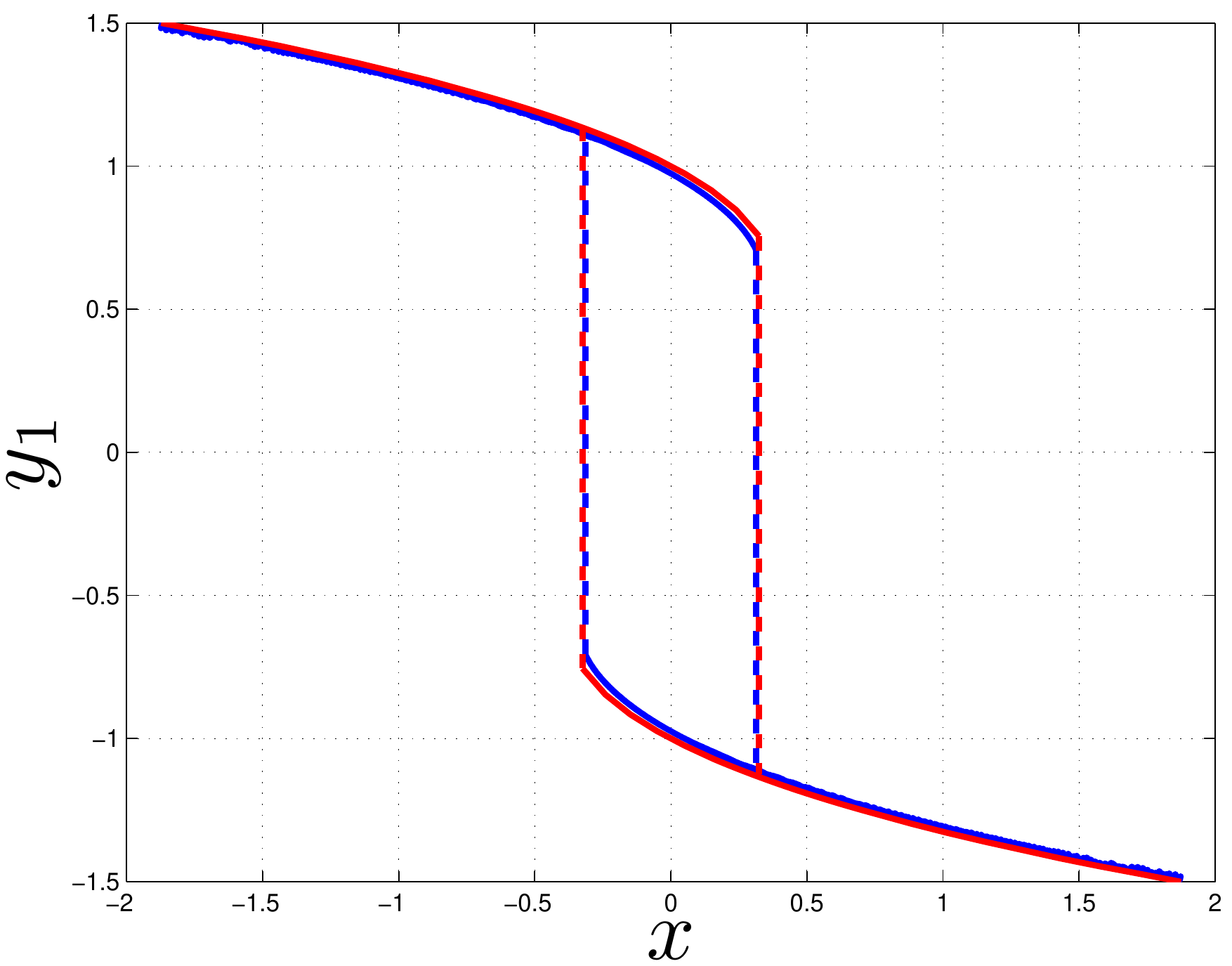}
\makeatletter
\def\@captype{figure}
\makeatother
\caption{\small Comparison of time averaged dynamics with PTA dynamics. Blue curve indicates $\overline{y}_1$ vs. $x$; red curve indicates $\hat{y}_1$ vs. $x$. Solid line indicates slow dynamics; dash line indicates fast dynamics.} \label{Invariant_Arstein52}
\end{figure}

Figures~\ref{Invariant_Arstein52} and~\ref{Invariant_Lorenz} present the time averaged variables from running the slow-fast system and the limit slow variables, of the model problems in Section 3.1.1 and 3.1.2, respectively. For both cases, we observe good agreement between the time averaged dynamics and the averaged limit dynamics. It is worth mentioning that:
\begin{enumerate}
\item For the `Forced' Lorenz system, we successfully approximate the coarse dynamics starting from `arbitrary' fine initial conditions. The coarse dynamics is insensitive to the fine initial conditions, namely, starting from a point that is far away from the $\omega$-limit set, the fast dynamics will jump to the set quickly and evolve around it.
\item We also test $\hat{x^2}$ as a nonlinear function of the fast variables, where the result is consistent with $\overline{x^2}$, as shown in Figure~\ref{Invariant_Lorenz}(d); In contrast, using the Tikhonov approach and defining $x^2$ as a coarse variable of interest, we will finally obtain the limit dynamics of $\hat{x}^2$ rather than $\hat{x^2}$.
\item As $\varepsilon\rightarrow0$, the fast dynamics of the problem~(\ref{arstein52}) converges to periodic limit cycles, and the fast dynamics of the `Forced' Lorenz system exhibits chaotic patterns; However, in both cases, the time averaged coarse variables approximate the time average of the limit solutions (which retain fast motions).
\item For the problem in Section 3.1.1, the limit cycles in the upper and lower branch center around $z_{up}(x)$ and $z_{lo}(x)$, respectively, where the centers satisfy $H=0$ and the time averaged trajectory and the DAE dynamics follows this slow manifold of the line of centers; the averaged trajectory of the `Forced' Lorenz system is distinct for the three initial conditions that satisfy fast equilibrium, which causes the mismatch between the time averaged dynamics and the DAE dynamics starting from $H=0$.
\end{enumerate}
\begin{center}
\includegraphics [width=.8\columnwidth]{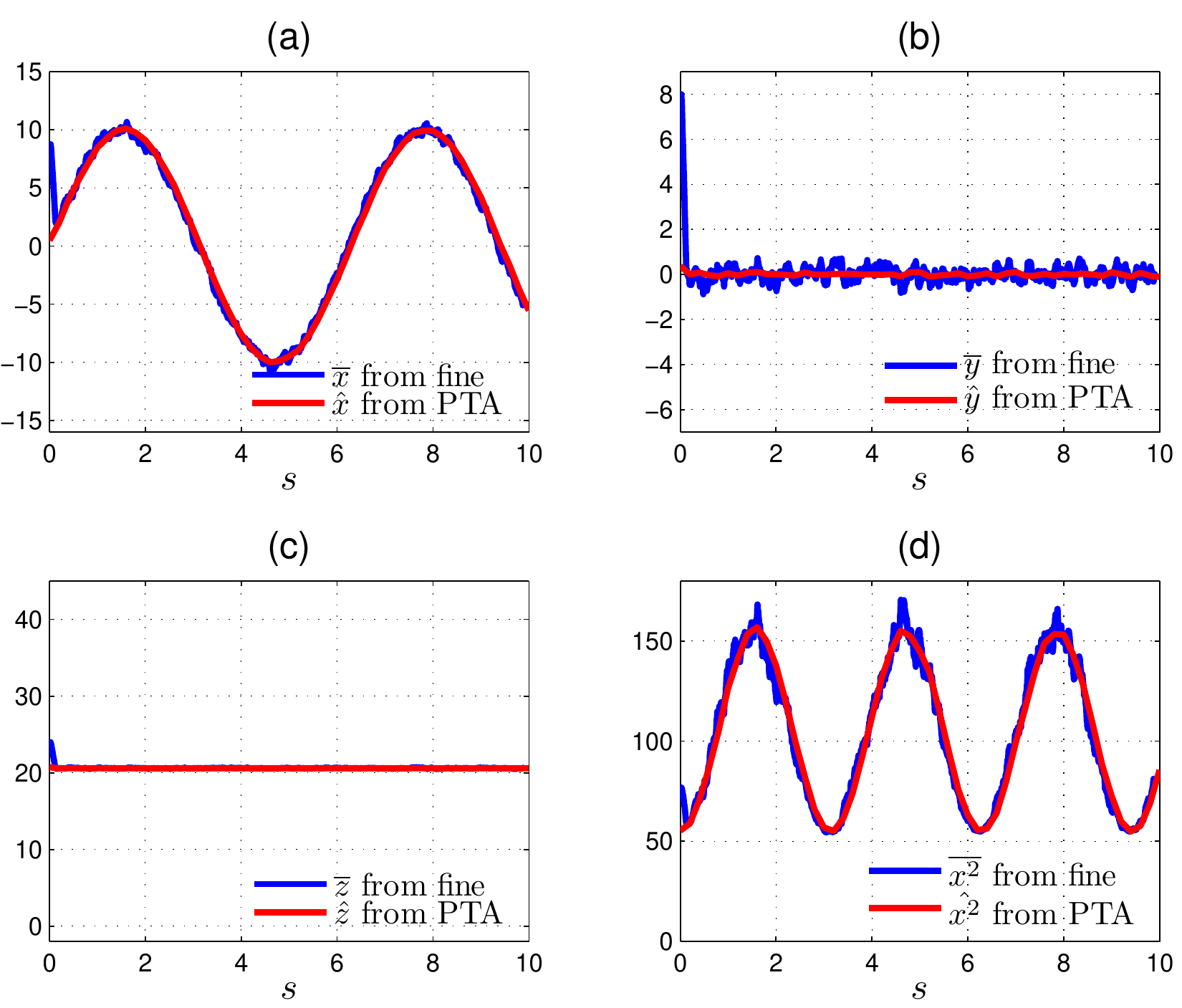}
\makeatletter
\def\@captype{figure}
\makeatother
\caption{\small Comparison of time averaged (coarse) variables from fine and from Practical Time Averaging (PTA)} \label{Invariant_Lorenz}
\end{center}
To show the necessity of defining `orthogonal observables' as coarse variables for the Young measure theory to be applicable, we follow~(\ref{instantaneous}) to write down the coarse evolution equations for instantaneous functions that are generally not orthogonal. We test this idea on `Forced' Lorenz system with the chosen coarse variables being $x$, $y$, $z$ and $x^2$. The corresponding coarse evolution equations are
\begin{equation}\label{non_orth}
\begin{aligned}
&\frac{dx}{ds}=\frac{1}{\omega}(\sigma(y-x)+\sigma{L})+g \\
&\frac{dy}{ds}=\frac{1}{\omega}(\gamma{x}-y-xz-\gamma{L}+Lz) \\
&\frac{dz}{ds}=\frac{1}{\omega}(xy-\beta{z}-Ly) \\
&\frac{dx^2}{ds}=2x\left(\frac{1}{\omega}(\sigma(y-x)+\sigma{L})+g \right)\\
&\frac{dL}{ds}=g\quad\quad\frac{dg}{ds}=-L.
\end{aligned}
\end{equation}
The results from running the system of ODEs~(\ref{non_orth}) are shown in Figure~\ref{Instant_Lorenz}. The instantaneous functions $x$, $y$, $z$ and $x^2$ are quite oscillatory due to $1/\omega$ in the evolution equation. This high oscillation is inappropriate for a physical coarse variable whose evolution should be amenable to easy approximation on the slow time scale. On the other hand, the time averaged variables $\overline{x}$, $\overline{y}$, $\overline{z}$ and $\overline{x^2}$ are much smoother and seem to be reasonable slow variables. This is one example of how PTA has the capability of producing robust slow variables for coarse dynamics.\\
\begin{center}
\includegraphics [width=.8\columnwidth]{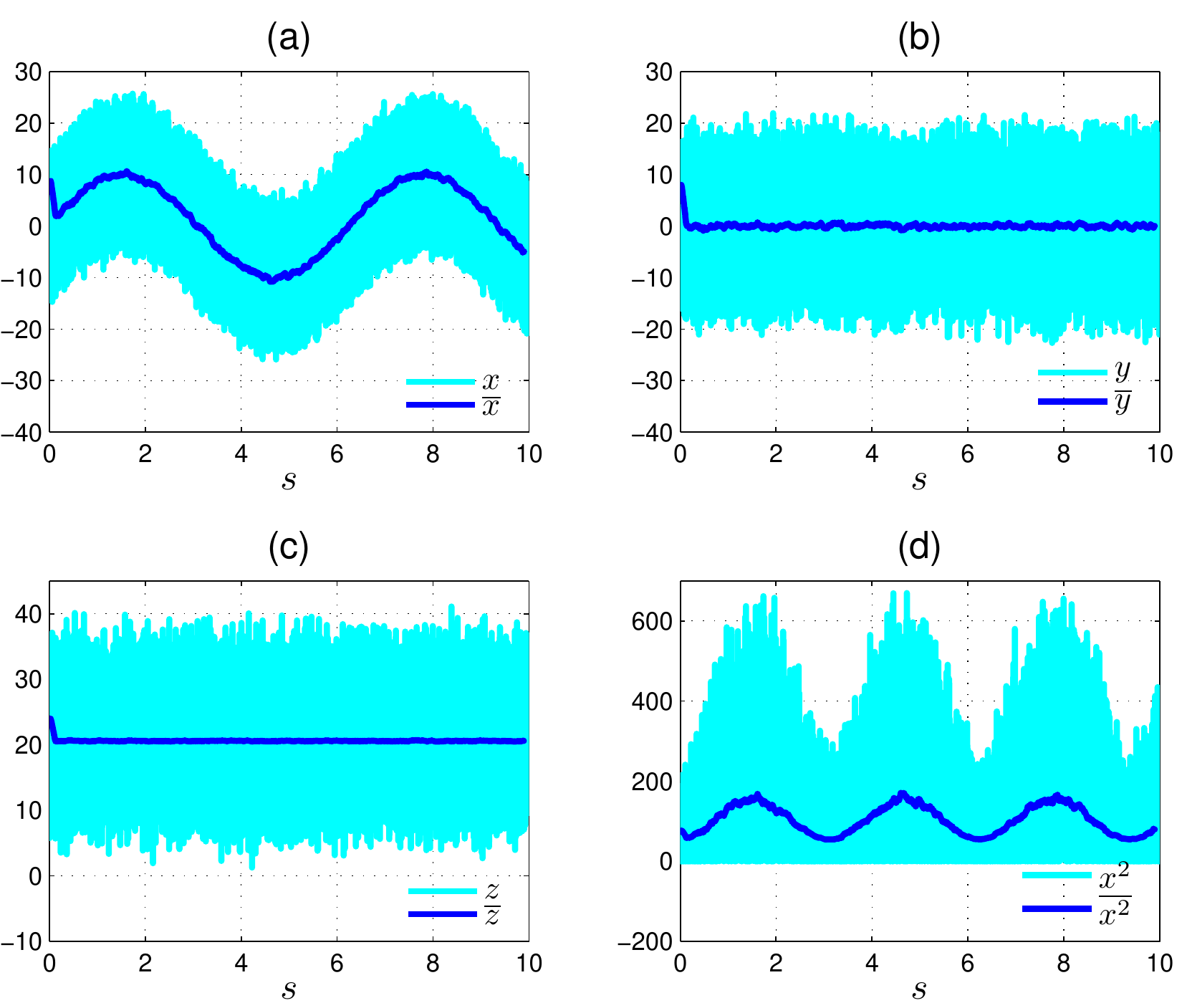}
\makeatletter
\def\@captype{figure}
\makeatother
\caption{\small Comparison of non-orthogonal instantaneous variables and time averaged variables on the slow time scale} \label{Instant_Lorenz}
\end{center}

\noindent
{\bfseries 4. Variables for autonomous coarse response}\\

\noindent
In the construction of coarse theory for fine systems, autonomous coarse response is usually desirable. However, the coarse theory based on the idea of PLIM or Practical Time Averaging cannot be guaranteed to be autonomous, i.e., the right-hand side term in the coarse theory is not uniquely defined from one single coarse state. The objective of this part is to explore this question and mathematically formulate a PDE for a class of coarse variables to formally allow an autonomous coarse response, where we further explore some appealing features of the selected coarse variables.\\

\noindent
4.1 MATHEMATICAL SCHEME \\

\noindent
Suppose the problem of interest is
\begin{equation}
\begin{split}
&\frac{df}{dt}(t)=H(f(t)) \\
&f(0)=f_*,
\end{split}
\end{equation}
where $f\in{\Bbb{R}}^N$. Based on the work in~\cite{acharya2007choice}, we define the coarse variable $c(t)$ as an instantaneous scalar function on the fine phase space and denote it as $\Pi$. Thus we have
\begin{equation}
c(t):=\Pi(f(t))\Longrightarrow\dot{c}(t)=\frac{\partial\Pi}{\partial f}(f(t))H(f(t)).
\label{eq2}
\end{equation}
From previous work~\cite{sawant2006model} on constructing the coarse dynamics of the Lorenz system using PLIM, we realize that for an arbitrary coarse function $\Pi$, the coarse dynamics cannot be ensured to be autonomous. Figure~\ref{Coarse_Figure1} shows the fine-to-coarse conversion of the Lorenz system by arbitrarily choosing $x$ and $z$ as coarse variables. Figure~\ref{Coarse_Figure1}(a) displays the fine trajectory on $(x,y,z)$ space. The coarse trajectory is simply the projection of the fine trajectory on $x$-$z$ plane, where we observe self-intersections of the trajectory, which simply imply a history-dependent behavior of the coarse evolution.
\begin{center}
\includegraphics [width=.6\columnwidth]{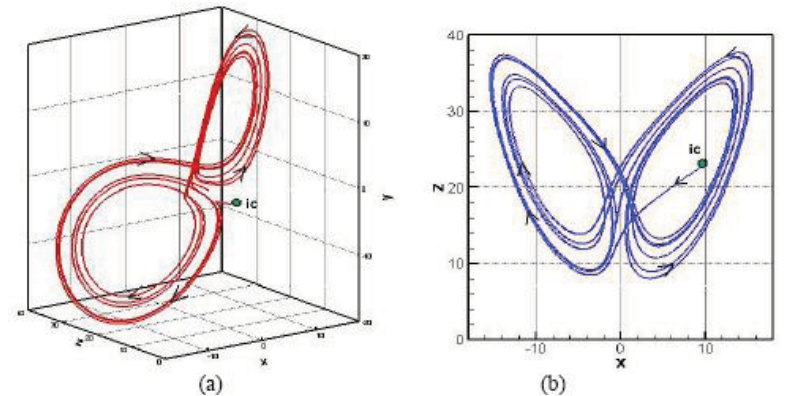}
\makeatletter
\def\@captype{figure}
\makeatother
\caption{\small Fine to coarse conversion (a) Fine response; (b) Coarse response}\label{Coarse_Figure1}
\end{center}
Therefore, the coarse variable needs to be specifically designed to formally ensure an autonomous coarse response. As outlined in~\cite{acharya2007choice}, a coarse variable is sought such that, if two fine states correspond to an identical coarse value, so does the coarse evolution from these fine states. Mathematically, this requires the existence of a function $\Pi$ for which there is a non-empty set of fine states defined by (\ref{eq3}), where $c$ is a constant and $A_c$ is another constant dependent on $c$ (and independent of fine states), i.e.
\begin{equation}
\{f:\Pi(f)=c\;\;\; \mbox{and}\;\;\; \frac{\partial\Pi}{\partial f}(f)H(f)=A_c\}.
\label{eq3}
\end{equation}
The existence of such a function (with enough smoothness) would imply autonomous response by (\ref{eq2}). Geometrically (\ref{eq3}) implies that a level set of $\Pi(\cdot)$ is also a level set of ${\partial\Pi}/{\partial f}(\cdot)H(\cdot)$. Therefore
\begin{equation}
\lambda\frac{\partial}{\partial f}(\Pi)=\frac{\partial}{\partial f}\left(\frac{\partial\Pi}{\partial f}H\right),
\label{eq4}
\end{equation}
where $\lambda$ is an arbitrary scalar. By taking $\lambda$ as a derivative of an arbitrary function $\Phi(\Pi)$ of the coarse variable, we have
\begin{equation}
\lambda:=\frac{\partial\Phi(\Pi)}{\partial\Pi}\Longrightarrow\frac{\partial}{\partial f}\left(\Phi(\Pi)-\frac{\partial\Pi}{\partial f}H\right)=0.
\label{eq5}
\end{equation}
One sufficient condition of~(\ref{eq5}) is to make $\Pi$ satisfy the first order PDE
\begin{equation}
\Phi(\Pi)=\frac{\partial\Pi}{\partial f}H.
\label{eq6}
\end{equation}
Thus, a function $\Pi$ satisfying~(\ref{eq6}) satisfies the necessary condition for autonomous response of $\Pi$, i.e.~(\ref{eq4}). In order to test this idea, we choose $\Phi$ to be an affine function of $\Pi$, namely,
\begin{equation}
\frac{\partial\Pi}{\partial f}H=a+b\Pi,
\label{eq7}
\end{equation}
where $a$ and $b$ are constants. (\ref{eq7}) is solved by implementing the finite element scheme. First, the coarse variable is discretized on the fine phase space, which takes the form
\begin{equation}
\Pi(f)=\sum_{A=1}^{N} \Pi^A\varphi^A(f).
\label{eq8}
\end{equation}
The residual of~(\ref{eq7}) is minimized, which leads to the following linear equations
\begin{equation}
\sum_{B=1}^{N} k_{AB}*\Pi^B=b_A,
\label{eq9}
\end{equation}
where
 \[k_{AB}=\int\left(\frac{\partial\varphi^A}{\partial f}H-b\varphi^A\right)\cdot\left(\frac{\partial\varphi^B}{\partial f}H-b\varphi^B\right)df \]
\[b_A=\int a\left(\frac{\partial\varphi^A}{\partial f}H-b\varphi^A\right)df. \]
We denote the linear equations in the matrix form
\begin{equation}\label{eq:linear_SVD}
Kx=b
\end{equation}
with
$K=\left(
\begin{array}{cccc}
k_{11} & k_{12} & \ldots & k_{1n}\\
k_{21} & k_{22} & \ldots & k_{2n}\\
\vdots & \vdots & \ddots & \vdots\\
k_{n1} & k_{n2} & \ldots & k_{nn}
\end{array} \right)$,
\quad $x=\left(
\begin{array}{c}
\Pi^{1} \\
\Pi^{2} \\
\vdots \\
\Pi^{N}
\end{array} \right)$\quad
and\quad $B=\left(
\begin{array}{c}
b_{1} \\
b_{2} \\
\vdots \\
b_{N}
\end{array} \right)$. \\
The system of linear equations is solved by the Singular Value Decomposition (SVD) method~\cite{press1992numerical} due to the singularity of the stiffness matrix. The stiffness matrix is decomposed as a product of three matrices, i.e.
\begin{equation}
K=U\cdot W\cdot V^T=U\cdot [diag(w_j)]\cdot V^T,
\end{equation}
where $U$ and $V$ are orthogonal, and $W$ is diagonal. For the singular matrix $K$, one or more $w_j's$  will be zero. If more $w_j's$  are zero, the matrix is more singular. The matrices computed from SVD simply give the orthonormal bases of the range and nullspace of the matrix $K$. Multiple solutions will be found because of the singularity of the stiffness matrix. Any solution to~(\ref{eq9}) defines an approximate solution to~(\ref{eq7}). We optimize over these solutions to satisfy the conditions mentioned below.\\

\noindent
4.2 PROBLEM DESCRIPTION \\

\noindent
Two dynamical systems are examined in our preliminary tests. One is the Lorenz system~\cite{lorenz1963deterministic}; the other is the Hald problem~\cite{chorin2000optimal}. \par
We recall the ODEs of the Lorenz system:
\begin{equation}
\begin{split}
&\dot{x}=\sigma(y-x) \\
&\dot{y}=\gamma x-y-xz \\
&\dot{z}=xy-\beta z \\
&\sigma=10, \;\beta=8/3,\; \gamma=25.
\end{split}
\end{equation}
The governing equation to solve for the coarse function takes the form
\begin{equation}
\frac{\partial\Pi}{\partial x}\cdot\sigma(y-x)+\frac{\partial\Pi}{\partial y}\cdot(\gamma x-y-xz)+ \frac{\partial\Pi}{\partial z}\cdot(xy-\beta z)=a+b\Pi.
\end{equation}
Figure~\ref{Coarse_Figure2} shows a level set of fine variables corresponding to an identical coarse value $c$. In this case, $c$ equals -0.0081. Also, we can see from the figure that there are multiple fine states corresponding to one single coarse state.
\begin{center}
\includegraphics [width=.5\columnwidth]{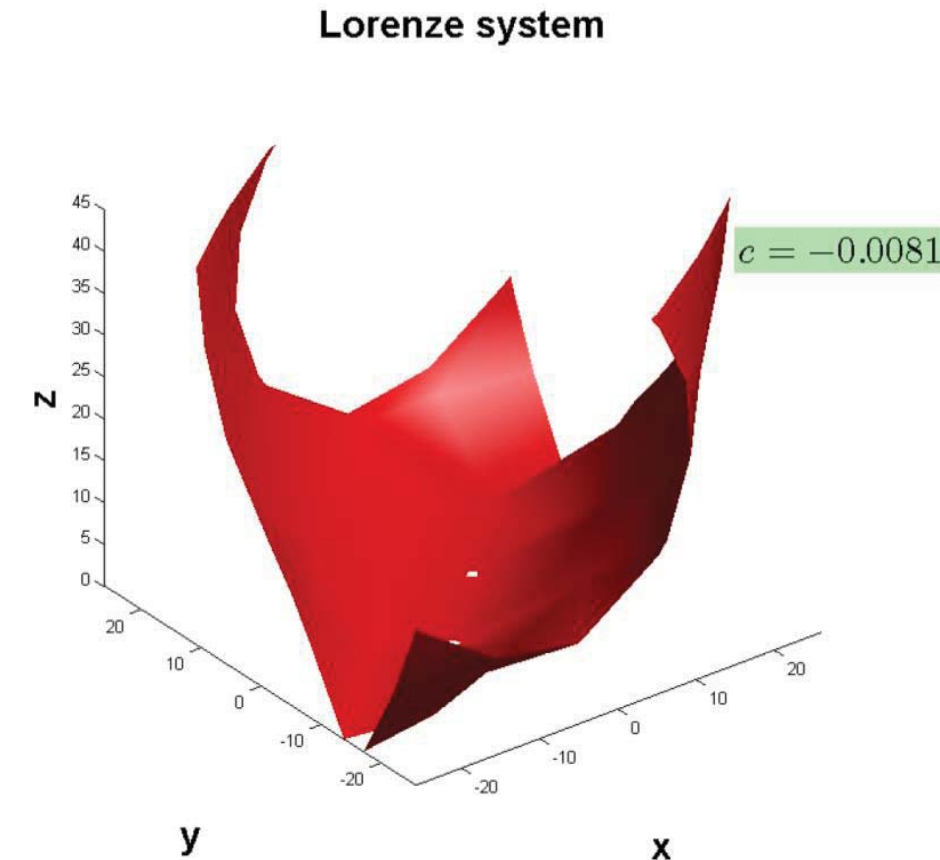}
\makeatletter
\def\@captype{figure}
\makeatother
\caption{\small A level set of fine variables}\label{Coarse_Figure2}
\end{center}
The second problem is a 4D Hamiltonian system proposed by Chorin et al~\cite{chorin2000optimal}. The ODEs of this system are
\begin{equation}
\begin{split}
&\dot{x_1}=x_2 \\
&\dot{x_2}=-x_1(1+x^2_3) \\
&\dot{x_3}=x_4 \\
&\dot{x_4}=-x_3(1+x^2_1),
\end{split}
\end{equation}
where $(x_1,x_2)$ and $(x_3,x_4)$ represent the displacements and velocities of two nonlinearly coupled oscillators.  Similarly, the PDE of the Hald problem is of the form
\begin{equation}
\frac{\partial\Pi}{\partial x_1}\cdot x_2+\frac{\partial\Pi}{\partial x_2}\cdot(-x_1(1+x^2_3))+\frac{\partial\Pi}{\partial x_3}\cdot x_4+
\frac{\partial\Pi}{\partial x_4}\cdot(-x_3(1+x^2_1))=a+b\Pi.
\end{equation}
The ultimate goal of choosing coarse variable is to generate autonomous coarse response. Notice that the mathematical statement we present is a {\it necessary} but not sufficient condition towards this goal. Moreover, one can hardly find an analytical solution to~(\ref{eq7}), and the numerical approximations do not have $C^1$ continuity. So the coarse variable needs to be carefully selected from the number of solutions to (\ref{eq9}) so that we obtain desirable characteristics in coarse response `close'\;to being autonomous. This selection is based on three criteria: \\
\noindent
Criterion (1). {\bfseries Consistency of coarse trajectories associated with symmetries of the fine system:} This criterion means the following: suppose the fine system is such that there exists a non-empty set $C$ of trajectories which is endowed with a non-empty (for simplicity, countable) set of mappings
\begin{equation}
Y_C = \{s_C: s_C(f^1(t)) = s_C(f^2(t))\; \forall\; t\quad\mbox{if}\;f^1\;\mbox{and}\;f^2\;\mbox{belong to}\;C  \}.
\end{equation}
It is possible that there could be more than one such set $C$ for a given fine system. We require that our computational procedure is successful in `blindly' capturing as many such functions $s_C$ as possible, i.e. without explicit knowledge of `symmetry' class $C$ built into the calculation of $s_C$.\\
\noindent
Criterion (2). {\bfseries Satisfactory match of coarse evolutions from distinct fine initial conditions with constraints:} Since the coarse variable is defined as a function of the fine variables, there exist multiple fine states associated with one coarse state and its time derivative. We deem a $\Pi$ as acceptable if for any two fine trajectories $f^1$ and $f^2$ that satisfy
\begin{equation}
\begin{aligned}
&\Pi\circ f^1(0)=\Pi\circ f^2(0)\\
&\frac{d}{dt}(\Pi\circ f^1)(0)=\frac{d}{dt}(\Pi\circ f^2)(0).
\end{aligned}
\label{constrain26}
\end{equation}
The difference of $\Pi\circ f^1$ and $\Pi\circ f^2$ is much smaller, in a suitable norm defined later, than the difference of $f^1$ and $f^2$.\\
\noindent
Criterion (3). {\bfseries Local convergence of coarse evolutions over a short time:} It is natural to require that $\Pi$ be such that whenever~(\ref{constrain26}) is satisfied, $\Pi\circ f^1$ and $\Pi\circ f^2$ agree well over short periods of time. \par

As discussed in the previous section, (\ref{eq:linear_SVD}) has numerous solutions since any vector in the nullspace can be added to a solution $x$. Therefore, the solution is a linear combination of two parts: the particular part and the homogeneous part, i.e.
\begin{equation}
x=x_{part}+\alpha x_{hom},
\end{equation}
where $\alpha$ is an arbitrary scalar. $x_{part}$ denotes any particular solution to~(\ref{eq:linear_SVD}), and $x_{hom}$ is any vector in the nullspace. The smoothness of the coarse function can be estimated by comparing the ratio of the variance to the mean of $x_{part}$ and $x_{hom}$, where we find that $x_{part}$ generally fluctuates more than $x_{hom}$. Thus, by changing the value of $\alpha$ we can `adjust' the smoothness of the coarse function. If $\alpha$ is larger, the coarse function is smoother. If $\alpha$ goes to infinity, the coarse function is close to a constant. We will not choose an infinitely large $\alpha$ since a constant coarse function is trivial. However, by modifying $\alpha$ we can obtain a coarse function with desirable smoothness so as to fulfill the requirement of Criteria (2) and (3). \par

A measure is introduced to estimate the difference between a pair of trajectories starting from distinct fine initial conditions, as illustrated in the equations below
\begin{subequations}
\begin{equation}
\triangle f:=\frac{\overline{\|f^1-f^2\|}}{1/2\cdot(\overline{\|f^1\|}+\overline{\|f^2\|})}
\end{equation}
\label{measure_a}
\begin{equation}
\triangle c:=\frac{\overline{\|c^1-c^2\|}}{1/2\cdot(\overline{\|c^1\|}+\overline{\|c^2\|})},
\end{equation}
\label{measure_b}
\end{subequations}
where
\begin{equation}
\|(.)\|:=\sqrt{\sum^N_{i=1}(.)^2_i}\quad \mbox{and} \quad \overline{(.)}:=\frac{1}{\mathbf{T}}\int^{\mathbf{T}}_{0}(.)(t)dt.
\end{equation}
$\mathbf{T}$ denotes the total simulation time. Here $f^1$ and $f^2$ refer to the fine trajectories starting from two distinct fine initial conditions, respectively. Similarly, $c^1$ and $c^2$ present the corresponding coarse trajectories, i.e. $c^i(t)=\Pi(f^i(t))$. This measure gives a sense of how much the trajectories differ from each other. The next section makes use of these measures to discuss the simulated results.\\

\noindent
4.3 NUMERICAL TESTS \\

\noindent
We test the properties of the computed coarse variables and present numerical results in this section.\\

\begin{figure}
\centering
\includegraphics [width=.6\columnwidth]{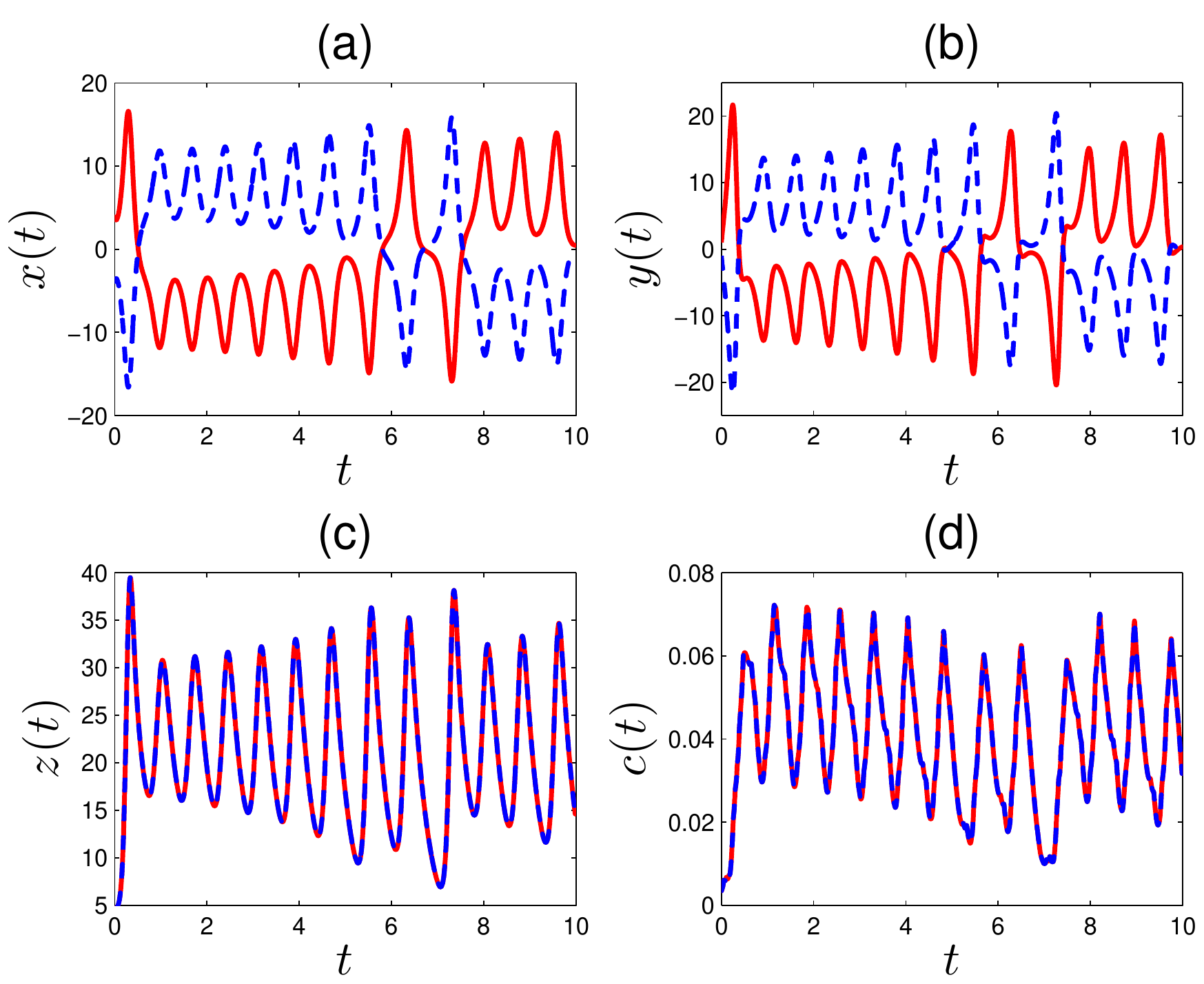}
\makeatletter
\def\@captype{figure}
\makeatother
\caption{\small Test symmetry on the Lorenz system (a) $x(t)$; (b) $y(t)$; (c) $z(t)$; (d) $c(t)$}\label{Coarse_Figure3}
\end{figure}

\noindent
(1). {\bfseries Consistency of coarse trajectories associated with symmetries of the fine system:} One symmetry of the Lorenz system is
\begin{equation}
\boxed{
\begin{array}{c}
x^2(0)=-x^1(0) \\
y^2(0)=-y^1(0) \\
z^2(0)=z^1(0)
\end{array}}
\longrightarrow
\boxed{
\begin{array}[c]{c}
x^2(t)=-x^1(t) \\
y^2(t)=-y^1(t) \\
z^2(t)=z^1(t).
\end{array}}
\label{eq15}
\end{equation}
Figures~\ref{Coarse_Figure3}(a),~\ref{Coarse_Figure3}(b) and~\ref{Coarse_Figure3}(c) show the fine evolutions. \emph{The trajectories with distinct fine initial values are indicated with different colors. The rule applies to all following plots.} For $x$ and $y$ variables, the trajectories are symmetric about $x=0$ and $y=0$, respectively, while the evolution of $z$ is identical in both cases. Figure~\ref{Coarse_Figure3}(d) shows the corresponding coarse evolutions, where the trajectories totally overlap each other. $c$ is the defined coarse variable $\Pi(f)$ as indicated in~(\ref{eq8}). It is remarkable that a discretely computed quantity preserves the symmetry property described in Criterion (1). Our procedure may be considered a generalization of the idea of finding invariants for autonomous ODE systems as, e.g., those considered in~\cite{artstein2009analysis}.\par
Suppose we choose a trajectory $(x^1(t),y^1(t),z^1(t))$. Then~(\ref{eq15}) and any one (or all) of the following functions $s_C(f):=|x|$, or $s_C(f):=|y|$, or $s_C:=\{|x|,\;|y|,\;{z}\}^T$ defines a class $C$ and the associated set of mappings $Y$ being discussed here. With these definitions in hand, we observe from Figure~\ref{Coarse_Figure3}(d) that $\Pi(f^1(t))=\Pi(f^2(t))\;\forall\; t$ for $f^1$ and $f^2$ belonging to $C$.\par
For the Hald problem, the symmetry has more than one form
\begin{subequations}
\begin{equation}\label{hald1}
\boxed{
\begin{array}{c}
x^2_1(0)=x^1_1(0) \\
x^2_2(0)=x^1_2(0) \\
x^2_3(0)=-x^1_3(0) \\
x^2_4(0)=-x^1_4(0)
\end{array}}
\longrightarrow
\boxed{
\begin{array}{c}
x^2_1(t)=x^1_1(t) \\
x^2_2(t)=x^1_2(t) \\
x^2_3(t)=-x^1_3(t) \\
x^2_4(t)=-x^1_4(t),
\end{array}}
\end{equation}
\begin{equation}\label{hald2}
\boxed{
\begin{array}{c}
x^3_1(0)=-x^1_1(0) \\
x^3_2(0)=-x^1_2(0) \\
x^3_3(0)=-x^1_3(0) \\
x^3_4(0)=-x^1_4(0)
\end{array}}
\longrightarrow
\boxed{
\begin{array}{c}
x^3_1(t)=-x^1_1(t) \\
x^3_2(t)=-x^1_2(t) \\
x^3_3(t)=-x^1_3(t) \\
x^3_4(t)=-x^1_4(t),
\end{array}}
\end{equation}
\begin{equation}\label{hald3}
\boxed{
\begin{array}{c}
x^4_1(0)=-x^1_1(0) \\
x^4_2(0)=-x^1_2(0) \\
x^4_3(0)=x^1_3(0) \\
x^4_4(0)=x^1_4(0)
\end{array}}
\longrightarrow
\boxed{
\begin{array}{c}
x^4_1(t)=-x^1_1(t) \\
x^4_2(t)=-x^1_2(t) \\
x^4_3(t)=x^1_3(t) \\
x^4_4(t)=x^1_4(t).
\end{array}}
\end{equation}
\end{subequations}
Once the fine initial values satisfy any of the possible symmetry forms, the fine evolutions will maintain the same kind of symmetry, as shown in Figure~\ref{Coarse_Figure4}. In this case, if we choose $s_C(f):=\{|x_1|,\;|x_2|\}^T$, or $s_C(f):=\{|x_3|,\;|x_4|\}^T$, or $s_C(f):=\{|x_1|,\;|x_2|,\;|x_3|,\;|x_4|\}^T$, these definitions of $s_C$ and classes of trajectories generated from any one of~(\ref{hald1}),~(\ref{hald2}),~(\ref{hald3}) serve to define autonomous coarse functions that we desire. \par
\begin{center}
\includegraphics [width=0.7\columnwidth]{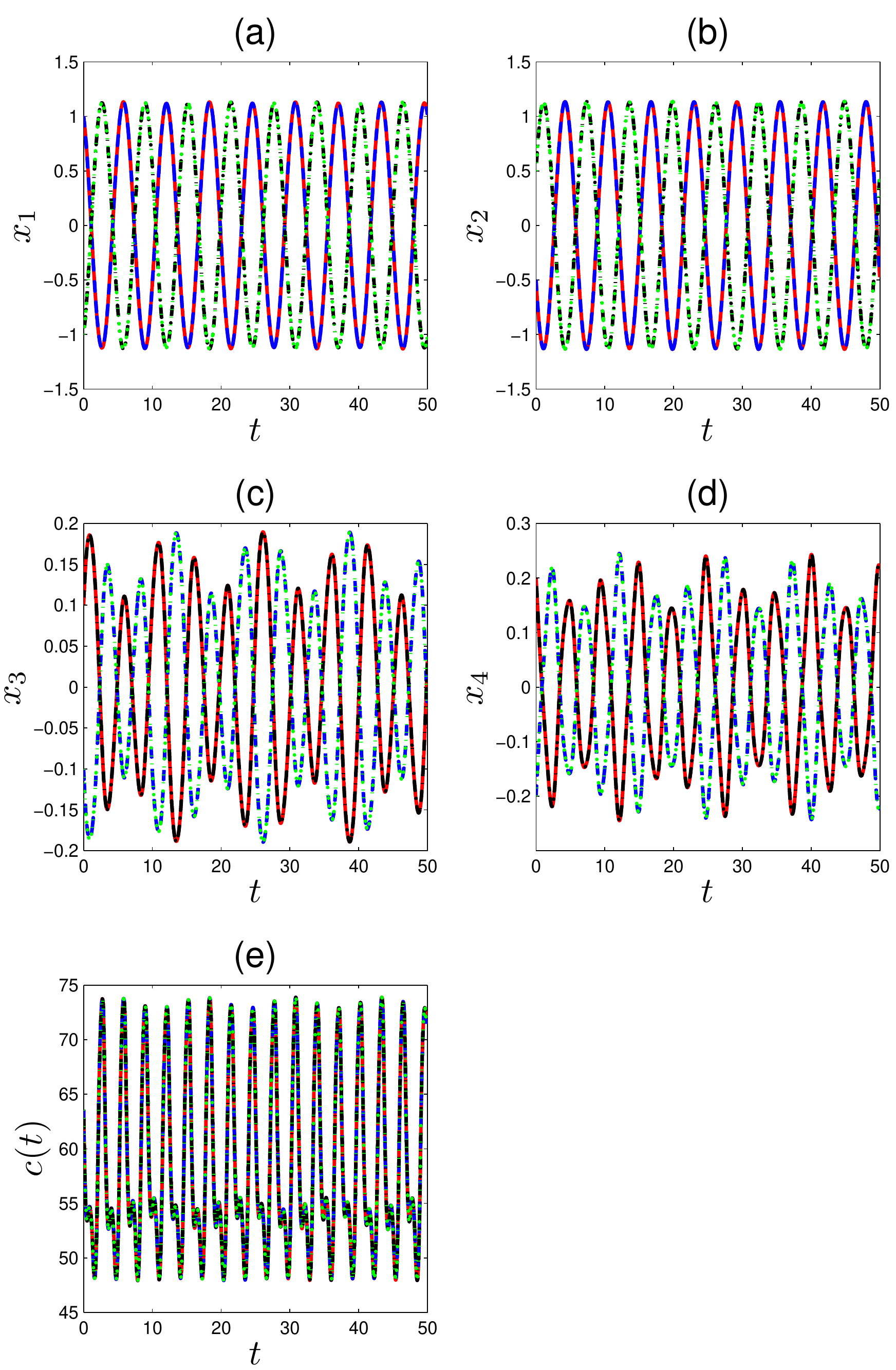}
\makeatletter
\def\@captype{figure}
\makeatother
\caption{\small Test symmetry on the Hald problem (a) $x_1(t)$; (b) $x_2(t)$; (c) $x_3(t)$; (d) $x_4(t)$; (e) $c(t)$}\label{Coarse_Figure4}
\end{center}
\noindent
(2). {\bfseries Coarse evolutions corresponding to distinct fine initial conditions with constraints:} We can see from Figure~\ref{Coarse_Figure5} that for fine trajectories that satisfy the constraints~(\ref{constrain26}), coarse evolutions match well, even though fine trajectories show a large divergence. By comparing the fine and coarse differences with the same measure, we observe a significantly larger fine difference. The ratio ${\triangle c}/{\triangle f}$ is equal to 0.1306/0.7449=0.1753, where ${\triangle f}$ and ${\triangle c}$ are defined in~(\ref{measure_a}a) and~(\ref{measure_b}b), respectively. The same constraints are imposed on the fine initial data on the Hald problem. Similar results are found as shown in Figure~\ref{Coarse_Figure6}, in which ${\triangle c}/{\triangle f}=0.2501/1.1310=0.2211$.\par

\begin{center}
\includegraphics [width=.7\columnwidth]{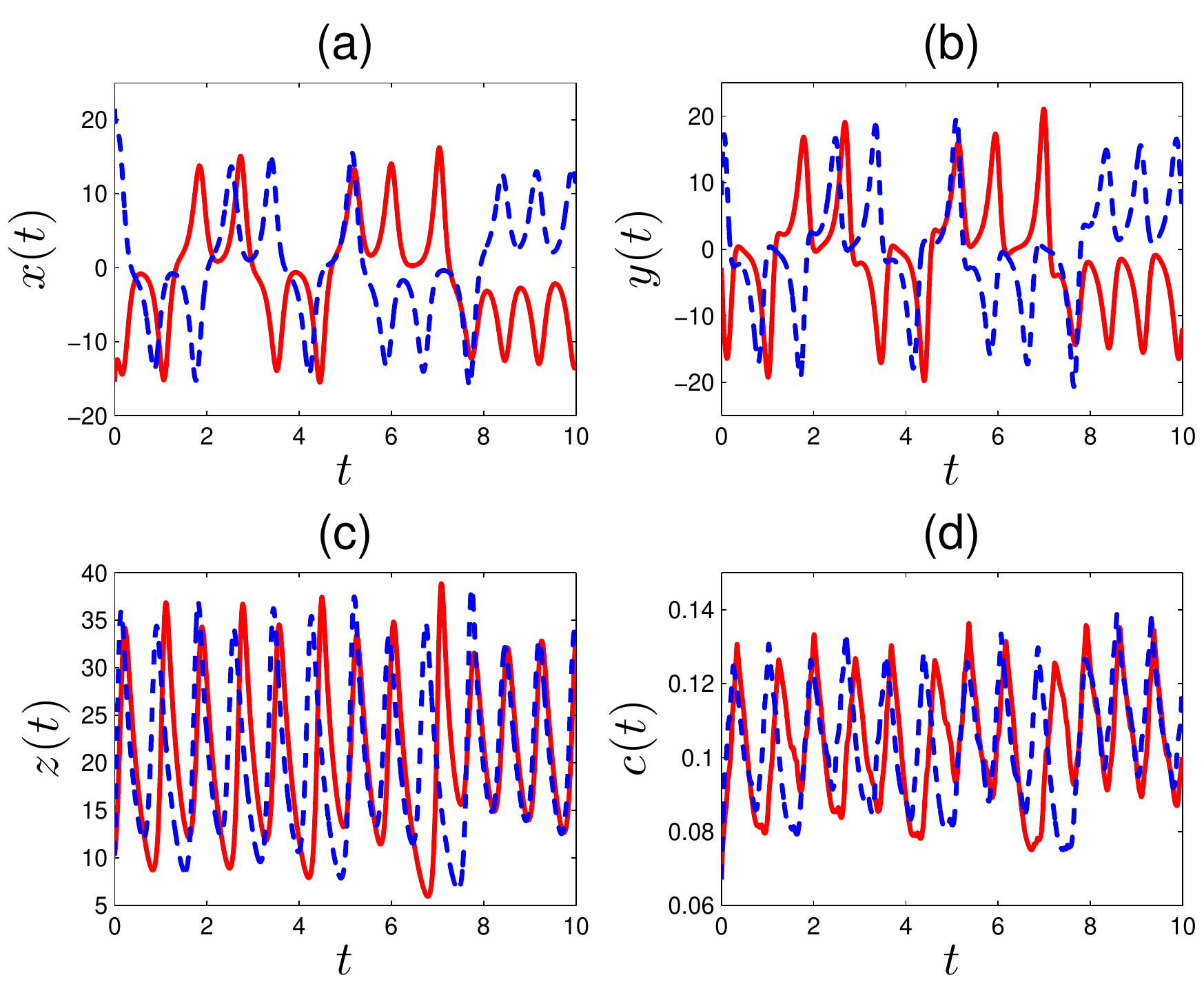}
\makeatletter
\def\@captype{figure}
\makeatother
\caption{\small Test coarse evolution with constraints on the Lorenz system (a) $x(t)$; (b) $y(t)$; (c) $z(t)$; (d) $c(t)$}\label{Coarse_Figure5}
\end{center}
\begin{center}
\includegraphics [width=.7\columnwidth]{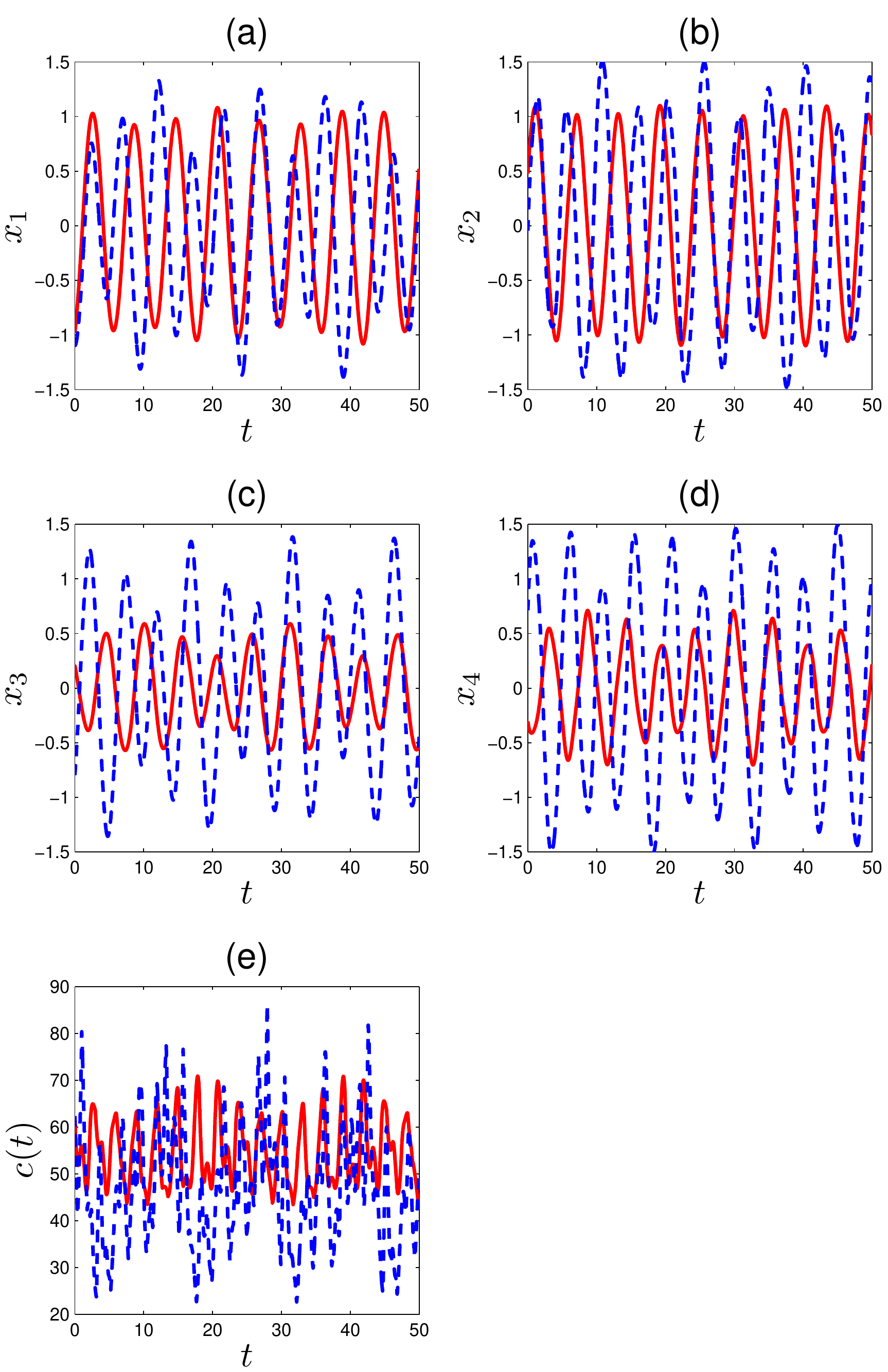}
\makeatletter
\def\@captype{figure}
\makeatother
\caption{\small Test coarse evolution with constraints on the Hald problem (a) $x_1(t)$; (b) $x_2(t)$; (c) $x_3(t)$; (d) $x_4(t)$; (e) $c(t)$}\label{Coarse_Figure6}
\end{center}
\noindent
(3). {\bfseries Local convergence of coarse evolutions over a short time:} Figure~\ref{Coarse_Figure7} shows four examples of the coarse evolutions within a shorter time. The constraints~(\ref{constrain26}) are imposed as before. We compute the coarse difference till $t=0.3$. The coarse differences $\triangle{c}$ are 0.0163, 0.0154, 0.1030 and 0.0765, respectively, which are generally smaller than the coarse differences for longer period of time shown in Figure~\ref{Coarse_Figure5}(d). However, for any two arbitrary fine initial conditions without these constraints, the coarse trajectories do not show local convergence. The coarse differences equal 0.2136, 0.2184, 0.2325 and 0.2224, respectively, which are more than 10 times bigger, as indicated in Figure~\ref{Coarse_Figure8} - this is simply to verify that this result is not simply an artifact induced by choosing a coarse variable that is a function on the fine phase space with negligible variation about a constant value. Similarly, Figure~\ref{Coarse_Figure9} shows the results of the Hald problem. We compute the coarse difference till $t=0.4$. For two coarse trajectories with the same initial value and its time derivative, the initial difference is considerably smaller compared with that from arbitrary fine initial conditions.\\
\begin{center}
\includegraphics [width=.7\columnwidth]{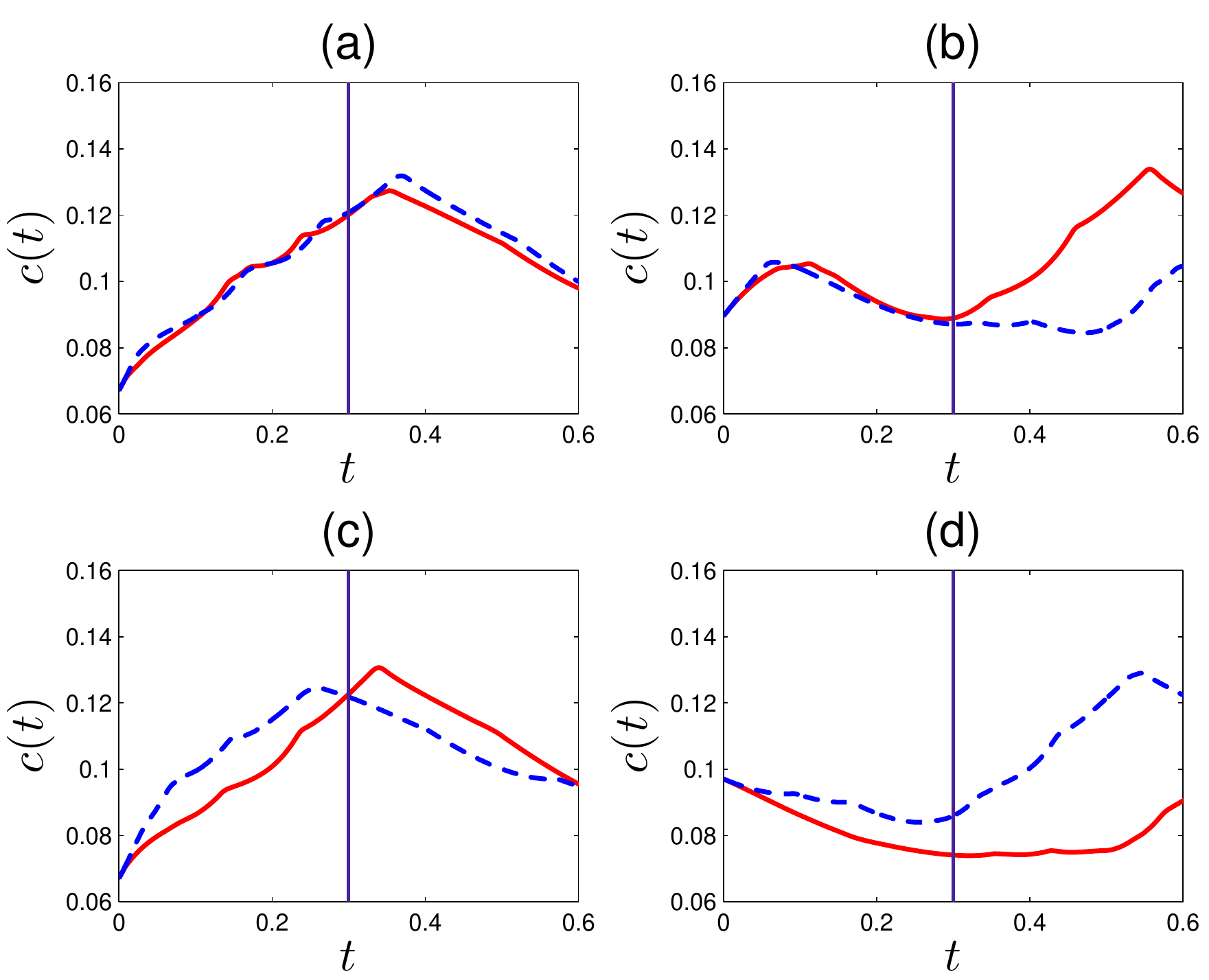}
\makeatletter
\def\@captype{figure}
\makeatother
\caption{\small Test local convergence on the Lorenz system (a) Case I: $\triangle c=0.0163$; (b) Case II: $\triangle c=0.0154$; (c) Case III: $\triangle c=0.1030$; (d) Case IV: $\triangle c=0.0765$}\label{Coarse_Figure7}
\end{center}
\begin{center}
\includegraphics [width=.7\columnwidth]{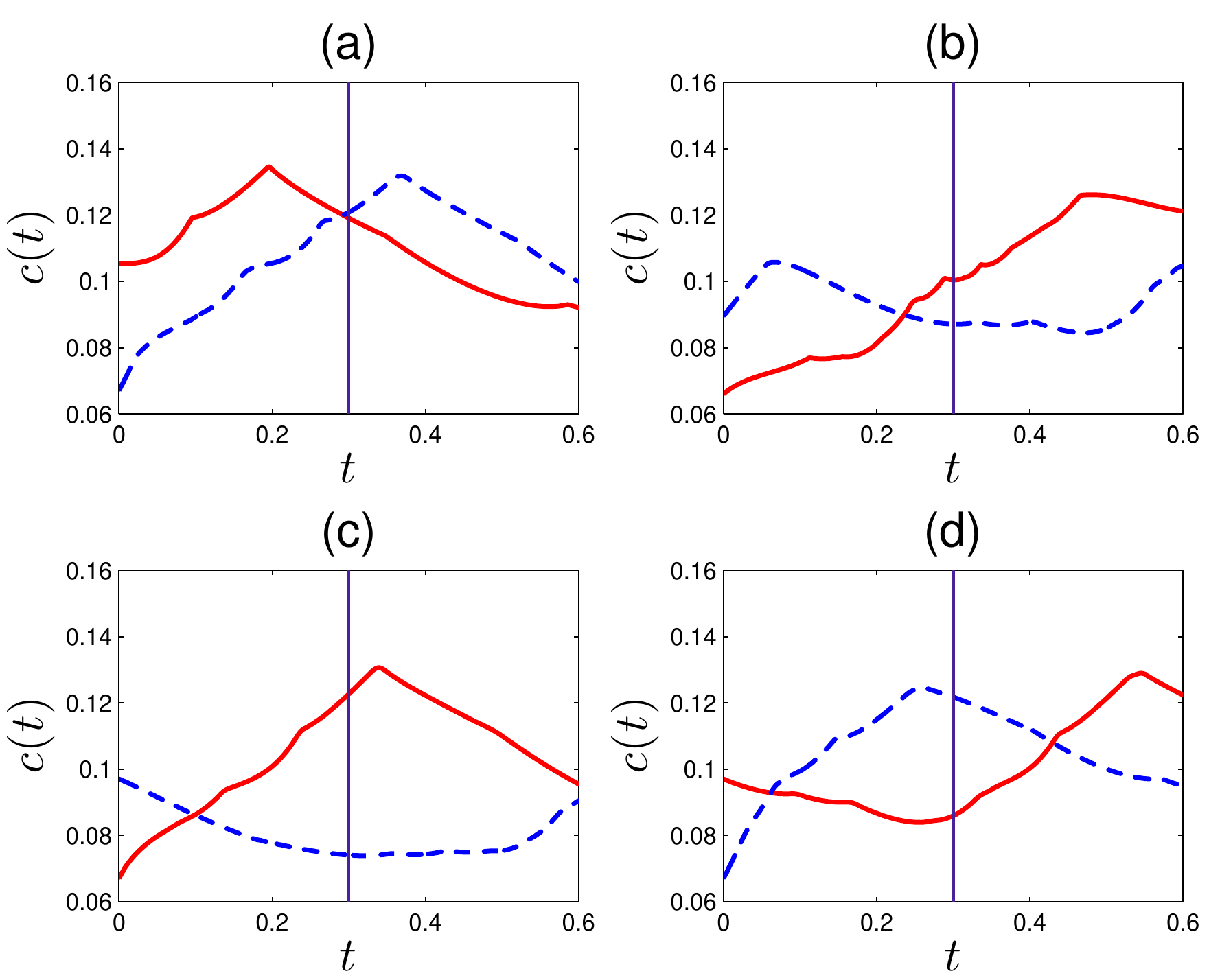}
\makeatletter
\def\@captype{figure}
\makeatother
\caption{\small $c(t)$ with arbitrary fine initial conditions on the Lorenz system (a) Case I: $\triangle c=0.2136$; (b) Case II: $\triangle c=0.2184$; (c) Case III: $\triangle c=0.2325$; (d) Case IV: $\triangle c=0.2224$}\label{Coarse_Figure8}
\end{center}

\begin{center}
\includegraphics [width=.7\columnwidth]{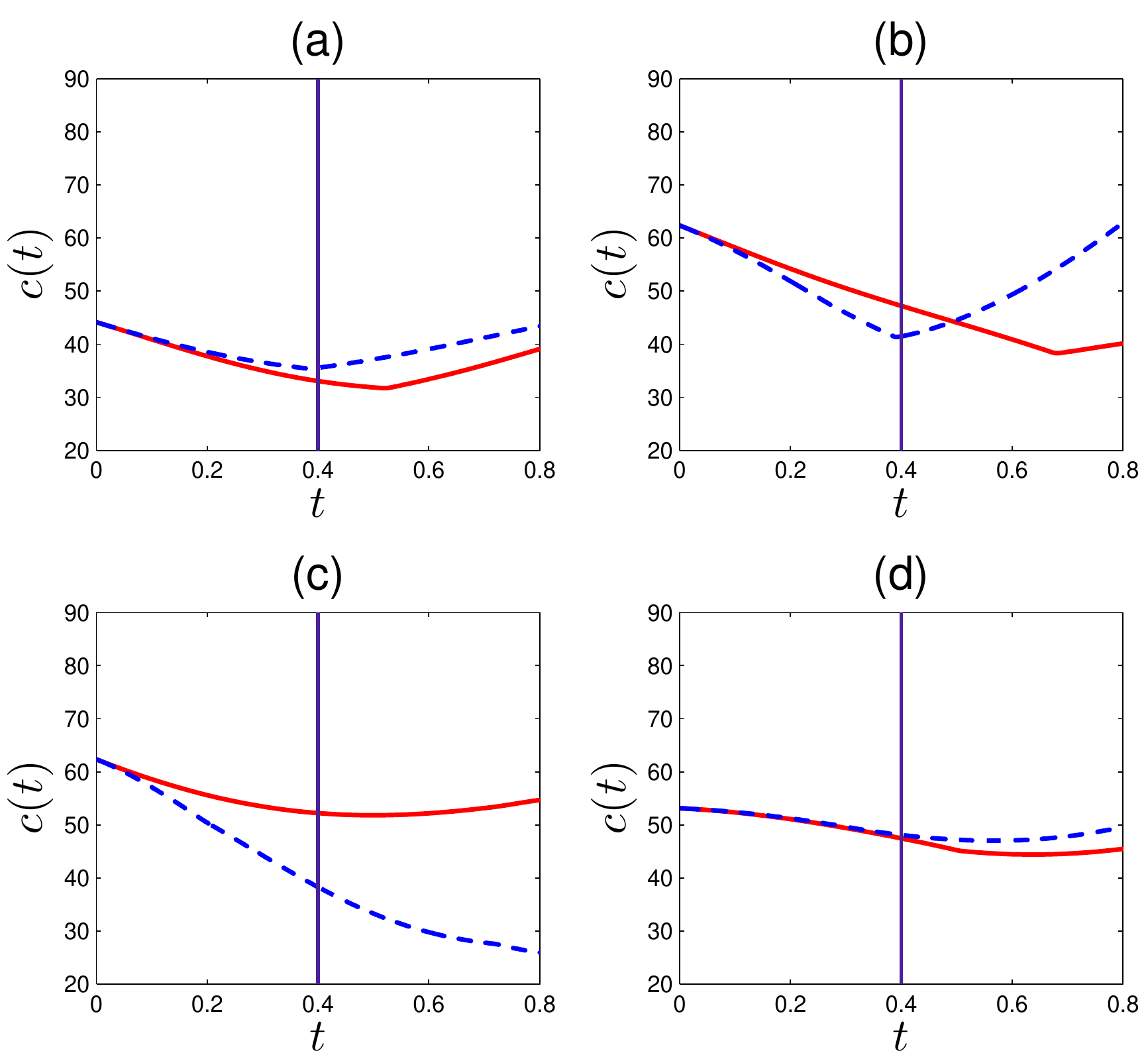}
\makeatletter
\def\@captype{figure}
\makeatother
\caption{\small Test local convergence on the Hald problem (a) Case I: $\triangle c=0.0231$; (b) Case II: $\triangle c=0.0498$; (c) Case III: $\triangle c=0.1042$; (d) Case IV: $\triangle c=0.0031$}\label{Coarse_Figure9}
\end{center}
\begin{center}
\includegraphics [width=.7\columnwidth]{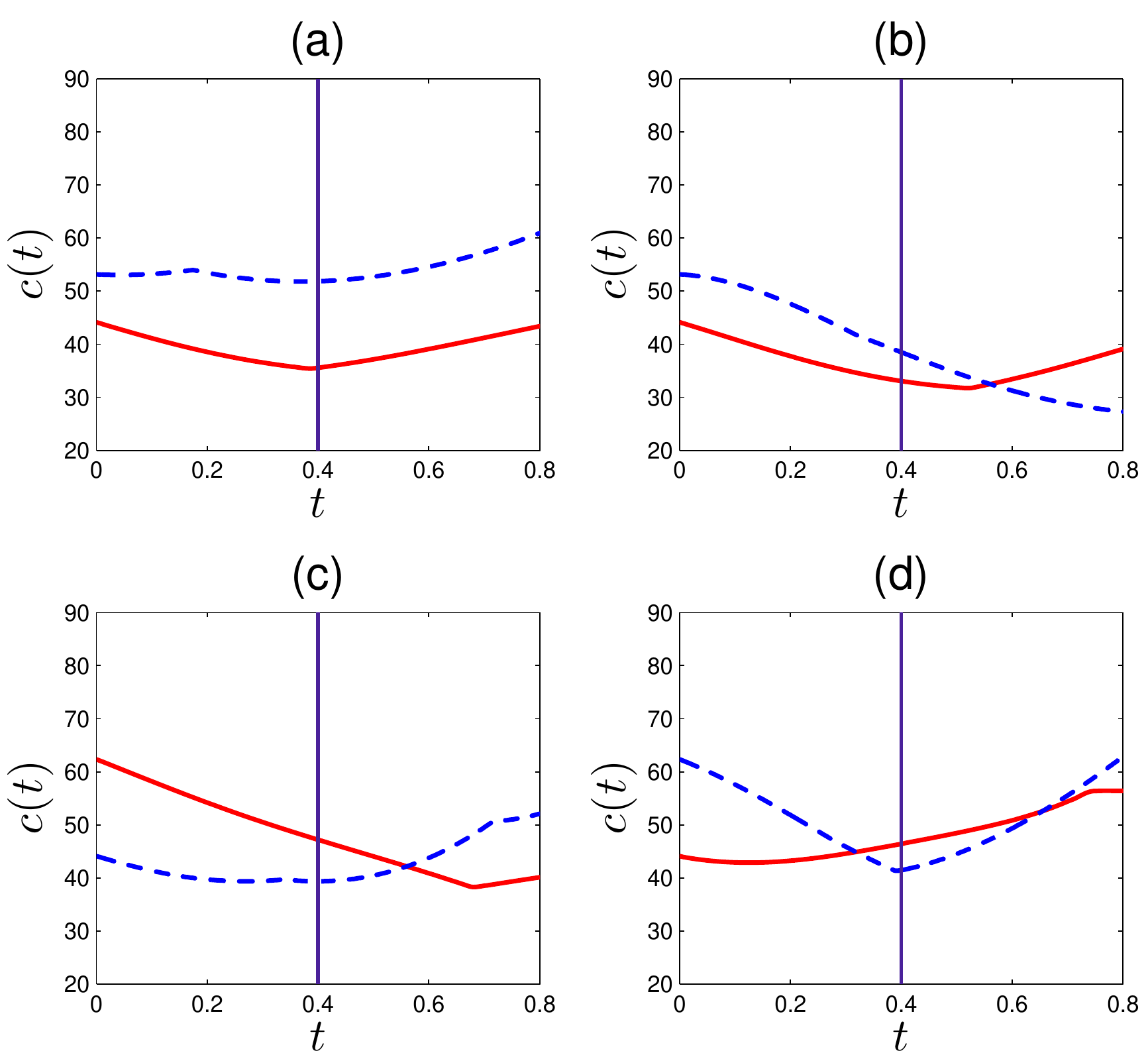}
\makeatletter
\def\@captype{figure}
\makeatother
\caption{\small $c(t)$ with arbitrary fine initial conditions on the Hald problem (a) Case I: $\triangle c=0.3019$; (b) Case II: $\triangle c=0.2097$; (c) Case III: $\triangle c=0.2929$; (d) Case IV: $\triangle c=0.1881$}\label{Coarse_Figure10}
\end{center}

\noindent
{\bfseries 5. Discussion}\\

\noindent
We summarize salient features of our work:
\begin{itemize}
\item Our main contribution is the sound demonstration of running-time averages as appropriate coarse variables of fast fine-scale dynamics. No other multi-scale computational method in the literature (``Equation-free'' approach~\cite{kevrekidis2003equation}, Heterogenous Multi-scale Method (HMM)~\cite{weinan2003heterogeneous, weinan}, Young measure theory~\cite{artstein2009analysis}) have discussed such a possibility. For example, the averaging of the discrete KDV-Burgers type equation in~\cite{artstein2009analysis} required the explicit knowledge of the $N/2+1$ invariants of the fast flow with period $N$; such information, unfortunately, is generally not available for dynamical systems.
\item We have explored three options for constructing coarse evolution equations. Each has its own pros and cons: \par
    \begin{itemize}
    \item In problems where time-averaged coarse response is insensitive to fine initial conditions, PLIM seems to perform reasonably. The timescale of the coarse dynamics constructed using PLIM is the averaging period $\tau$; while using the Tikhonov approach and Practical Time Averaging, we seek the average of the limit dynamics at the timescale $T_I$, which is the timescale of loading period and is usually provided by the physical problem. Suppose we have a MD system with atoms vibrating at timescale $t$ and an external loading applied to the system with loading period $T_I$, with measurements of different resolutions $\tau$. Normally we have $t<\tau<T_I$. PLIM allows us to `observe' the dynamics at different scales. On the other hand, in the study of the singularly perturbed system with $\varepsilon\rightarrow0$, the dynamics at the intermediate scale is not well defined. However, computation of the invariance PDE is a non-trivial task, especially when $\varepsilon\rightarrow0$. The computational methods we have tried are quite simplistic and it remains to be seen how Sacker's~\cite{sacker1965new} elliptic regularization existence guarantees fares for PLIM with a small number of coarse variables as well as more sophisticated numerical schemes tailored to hyperbolic PDE~\cite{guckenheimer2004fast} - the latter for the case when $\varepsilon$ is not far less than 1.
    \item For the limiting case when $\varepsilon \rightarrow 0$, the Tikhonov approach is efficient in computation, compared with solving the PDE in PLIM and provides a candidate ODE system (defined on the slow time scale with no $\varepsilon$) for PLIM to be further implemented. However, invertibility of the matrix ${\partial H}/{\partial f}$ in~(\ref{Approx_dynamics}) is a serious issue. It also has its limitations in that, at best, it can approximate time averaged dynamics robustly only when the fast dynamics settles on an equilibria. Moreover, even with latitude, it cannot be expected to approximate the time-averaged behavior of nonlinear functions of state (while, as demonstrated herein, it can on occasion do a reasonable job for linear state functions).
    \item Practical Time Averaging takes into account more general cases with the only requirement that the fast flow (for each fixed coarse state) remain in a bounded region of phase space. We have demonstrated several such cases in this paper. It has been shown that the limit behavior of nonlinear functions can be approximated as well. It is also efficient in computation since it allows a very big interval $T$ as the coarse time step.
    \end{itemize}
\item Construction of autonomous coarse system is still an open question. With the implementation of PLIM, the emergence of physically history-dependent behavior from underlying current-state-dependent behavior is observed (for instance, in~\cite{sawant2006model}), which arises as a direct consequence of reduction or coarse observation of the fine system. As discussed in detail in~\cite{slemrod2010time} such dependence is not expected to go away even when time-averaged coarse variables are used. For non-history dependent coarse variables, such non-autonomous coarse behavior was most explicitly noted in~\cite{artstein1996singularly} and subsequently in~\cite{artstein2007slow}. The computations in~\cite{artstein2009analysis} do not emphasize such dependence but, as pointed out in~\cite{slemrod2010time}, it is of course there and is accounted for by the special choice of initial conditions that is made to computationally generate the Young measure averages.
    This issue necessarily exists in HMM too in their reconstruction step. For instance, to obtain the stress (the only unknown data) in the macroscopic model of complex fluids (in which the coarse variables are the macroscopic velocity field and the stress), they~\cite{ren2005heterogeneous} apply the Irving-Kirkwood formula, which is dependent on the positions of particles. These quantities are defined in the microscopic model. From a physical point of view, fine quantities that correspond to a unique macroscopic velocity field can not be ensured to generate a unique stress.
\item In the second part of the paper, we consider the more classical question of coarse variables that are `instantaneous' functions of the fine variables. An interesting finding here is that we are able to \emph{numerically} define some such functions whose evolution is guaranteed to be unique when evaluated on whole classes of fine trajectories.
\end{itemize}

\noindent
{\bfseries Acknowledgement.} The work was partially supported by Civil and Environmental Engineering Department at Carnegie Mellon University through the Ellegood Strategic Doctoral Fellowship and partially supported by NSF (Dynamical Systems, Grant \# 0926579). \\

\bibliographystyle{plain}
\bibliography{ref}

\end{document}